\documentclass[11pt]{amsart}
\usepackage{amsmath,amssymb}
\numberwithin{equation}{section}

\usepackage[pagebackref,colorlinks=true,linkcolor=blue,citecolor=blue]{hyperref}

\usepackage{tikz}
\usetikzlibrary{arrows.meta}
\definecolor{Kblue}{RGB}{40,84,115}
\definecolor{Krust}{RGB}{151,77,50}
\definecolor{Kgreen}{RGB}{65,105,85}
\newtheorem{theorem}{Theorem}[section]
\newtheorem{lemma}[theorem]{Lemma}
\newtheorem{corollary}[theorem]{Corollary}

\newtheorem{definition}[theorem]{Definition}
\theoremstyle{remark}
\newtheorem{remark}[theorem]{Remark}
\newcommand{\R}{\mathbb R}
\newcommand{\tr}{\operatorname{tr}}

\newcommand{\diag}{\operatorname{diag}}
\newcommand{\cL}{\mathcal L}
\newcommand{\dd}{\,d}
\newcommand{\p}{\partial}
\newcommand{\Rn}{\mathbb R^n}
\newcommand{\la}{\lambda}
\newcommand{\diver}{\operatorname{div}}

\title[]
{Wiener's criterion for generalized Kolmogorov equations}
\date{}

\makeatletter
\@ifundefined{subjclassname@2020}{%
  \@namedef{subjclassname@2020}{\textup{2020} Mathematics Subject Classification}%
}{}
\makeatother


\author{Nicola Garofalo}

\address{School of Mathematical and Statistical Sciences\\ Arizona State University\\ Tempe, AZ 85287-1804}
\email{Nicola.Garofalo@asu.edu}

\begin{document}

\begin{abstract}
We establish a geometric criterion of Wiener type for generalized Kolmogorov equations
satisfying the
Kalman--H\"ormander condition. The
proof rests on a new strong Harnack inequality on kernel balls reminiscent of that proved by Evans and Gariepy in their breakthrough result for the heat equation. For degenerate operators the counterpart of the Evans--Gariepy strategy fails on a set exterior to any region where the ordinary Harnack inequality applies uniformly. We overcome this obstruction by a suitable correction, tailored to the geometry of the relevant operators.  
\end{abstract}

\maketitle

\tableofcontents

\section{Introduction and main theorem}\label{sec:introduction}

A basic problem in potential theory is the characterization of regular points, i.e., those boundary points at which the
generalized Perron--Wiener--Brelot solution continuously converges to its prescribed
Dirichlet datum.  
Wiener's classical geometric criterion \cite{Wie} characterized the
regular boundary points for the Laplacian in terms of the capacities of annular regions.
Fifty-eight years later, Evans and Gariepy \cite{EG} obtained a
corresponding geometric result for the heat equation. The long time interval between the two results reflects the substantial difficulties in passing from the stationary to the evolution setting. For instance, a famous result of Littman, Stampacchia and Weinberger \cite{LSW} states that the regular boundary points for any uniformly elliptic operator with bounded measurable coefficients coincide with those of the Laplacian. This remarkable invariance has no counterpart in the evolution setting. Petrovski\u{\i}'s classical example \cite{Pet} shows that a boundary point may be regular for $\p_t-\alpha\Delta$ and irregular for $\p_t-\beta\Delta$, with $\alpha\ne\beta$, $\alpha, \beta>0$.   

In the present work we establish, for a class of possibly highly degenerate evolution operators, a geometric Wiener test in the spirit of Evans and Gariepy. The test is formulated in terms of the level sets of the fundamental solution and provides a characterization of the regular boundary points for
the generalized Kolmogorov equation
\begin{equation}\label{operator-definition}
 \cL u=\partial_tu-\tr(Q\nabla_x^2u)
             -\langle Bx,\nabla_xu\rangle=0.
\end{equation}
Here $Q=Q^T\geq0$ and $B$ are real $n\times n$
matrices satisfying the Kalman--H\"ormander condition:
\begin{equation}\label{kalman}
\operatorname{Ker}Q\ \emph{ does not contain any nontrivial invariant subspace of}\ B^T.
\end{equation}
When $Q = I$ and $B=0$, we obtain as a special case the heat equation, but the condition \eqref{kalman} allows for a rich variety of partial differential equations of considerable theoretical and applied interest.
While for a discussion of these aspects we refer to the works \cite{BG, GTma, GTpa}, we emphasize that the equations in the class \eqref{operator-definition} can be quite degenerate.
One basic example is the kinetic
equation in $\R^3$, introduced by Kolmogorov in \cite{Kol}, 
\begin{equation}\label{kolmo}
\p_t u - \p_{xx}u -x\p_y u = 0,
\end{equation}
which is not parabolic in the classical sense but belongs to the class
\eqref{operator-definition}, with
\begin{equation}\label{mako}
Q = \begin{pmatrix} 1 & 0 \\ 0 & 0\end{pmatrix},\ \ \ \ B = \begin{pmatrix} 0 & 0 \\ 1 & 0\end{pmatrix},
\end{equation}
see also \eqref{badd} below.
The equation \eqref{kolmo} was introduced by Kolmogorov in connection with the kinetic theory of gases. Diffusion acts only in the velocity variable $x$, while the transport term transfers its regularizing effect to the position variable $y$. In his note Kolmogorov found the following smooth explicit fundamental solution of \eqref{kolmo} with pole at the origin
\begin{equation}\label{kolp}
p(x,y,t) = \frac{\sqrt3}{2\pi}\,t^{-2}\,
\exp\left[-\left(\frac{x^2}{4t}
+\frac{3}{t^{3}}\Big(y+\frac{t}{2}x\Big)^{2}\right)\right].
\end{equation}

It was H\"ormander who, in the opening of his 1967 work \cite{Hor}, placed Kolmogorov's fundamental solution \eqref{kolp} in
its general setting. After noting that for the class
\eqref{operator-definition}, condition \eqref{kalman} is equivalent to the
strict positivity of the Gramian matrix
\begin{equation}\label{Qt}
Q(t)=\int_0^t e^{sB}Qe^{sB^T}\dd s\ > 0,\ \ \ \ \ \ t>0,
\end{equation}
H\"ormander showed that, with the normalization
\begin{equation}\label{intrinsic-distance-volume}
\begin{cases}
m_t(x,y)^2 =\frac14\left\langle Q(t)^{-1}(e^{tB}x-y),e^{tB}x-y\right\rangle,\ \ \  t>0,
\\
V(t)=\sqrt{\det Q(t)},
\end{cases}
\end{equation}
the fundamental solution of \eqref{operator-definition}, with pole at
$(y,0)$, is given by
\begin{equation}\label{transition-kernel}
 p(x,y,t)=\frac{(4\pi)^{-n/2}}{V(t)}e^{-m_t(x,y)^2}.
\end{equation}
For the matrices in \eqref{mako}, formula \eqref{transition-kernel}
reduces to \eqref{kolp}, see  Section \ref{sec:model-geometry}.

It is worth noting here that the differential operator $\cL$ on the left-hand side of \eqref{operator-definition}  is  invariant under the non-commutative group law 
\begin{equation}\label{Lie}
(x,s)\circ (y,t) = (y+ e^{-tB}x,s+t),
\end{equation}
where $x, y\in \Rn$ and $s, t\in \R$, and such invariance is reflected in the time-varying pseudometric $m_t(x,y)$ in \eqref{intrinsic-distance-volume}. Equipped with \eqref{Lie}, the space $(\R^{n+1},\circ)$ forms a Lie group, whose identity element is the origin. Such group is non-Abelian when $B\neq0$.

Returning to the main objective of this paper, we mention that a geometric characterization of regular boundary points for
degenerate equations such as \eqref{kolmo}, and more generally
\eqref{operator-definition}, has remained open for quite some time. The
breakthrough result of Evans and Gariepy \cite{EG} for the heat equation,
and its extension in \cite{GL88} to uniformly parabolic equations with
variable coefficients, showed that the decisive ingredient is a strong
form of the Harnack inequality. A corresponding strong form of Harnack inequality has not been available
for the degenerate equations \eqref{operator-definition}.

The purpose of the present work is to fill this gap and use the resulting estimate to establish
the following geometric Wiener criterion.
To state our main results we need to introduce some notation.
 With $p$ as in \eqref{transition-kernel}, and $z=(x,t)$ and $\zeta = (y,\tau)$, we will write
 \begin{equation}\label{spacetime-kernel}
 \Gamma(z;\zeta)=
 \begin{cases}
 p(x,y,t-\tau),\ t>\tau,
 \\
 0,\ \ \ \ \ \ \ \ \ \ \ \ \ \ \ \ t\leq\tau,
 \end{cases}
 \end{equation}
to indicate the forward fundamental solution of \eqref{operator-definition} with pole at $\zeta$.
For $z_0=(x_0,t_0)$ and $r>0$, define the backward \emph{kernel ball} of
 radius $r$ and pole at $z_0$ by
\begin{equation}\label{kernel-ball-at-pole}
 \mathcal E_r(z_0)
 =\big\{\zeta\in\mathbb R^{n+1}\mid
 \Gamma(z_0;\zeta)>\frac{(4\pi)^{-n/2}}{V(r)}\big\}.
\end{equation}
Notice that \eqref{spacetime-kernel}
implies that $z_0\not\in \mathcal E_r(z_0)$. When $z_0 = (0,0)$, we simply write $\mathcal E_r$. The equations \eqref{intrinsic-distance-volume}-\eqref{transition-kernel} give the following equivalent description of \eqref{kernel-ball-at-pole}
\begin{equation}\label{kernel-ball-explicit}
 \mathcal E_r(z_0)
 =
 \big\{(y,\tau)\mid
 0<t_0-\tau<r,\ 
 m_{t_0-\tau}(x_0,y)^2
 <
 \log\frac{V(r)}{V(t_0-\tau)}\big\}.
\end{equation}

For each fixed time level $\tau\in(t_0-r,t_0)$, the spatial section of
$\mathcal E_r(z_0)$ is an ellipsoid centered at
$e^{(t_0-\tau)B}x_0$. The center tends to $x_0$ as
$\tau\to t_0^-$, and to $e^{rB}x_0$ as
$\tau\to(t_0-r)^+$. Its diameter tends to zero in both limits.
A distinct feature,
explained in Section~\ref{sec:scaling}, is that arbitrarily close to the
pole $z_0$ the kernel ball contains a portion which is exterior to any fixed region defined by the dilations \eqref{spacedil} attached to the canonical Kalman filtration. This accounts for the relevance of Theorem \ref{thm:strong} below.

Given a finite nonnegative Radon measure $\mu$, we define its potential by
\[
\Gamma_\mu(z)=\int_{\mathbb R^{n+1}}\Gamma(z;\zeta)\,d\mu(\zeta),
\]
where the integral is allowed to take the value $+\infty$. For a compact
set $F\subset\R^{n+1}$, define
\begin{equation}\label{capacity-definition}
 \operatorname{cap}(F)=\sup\bigl\{\mu(F)\mid \mu\geq0,
       \ \operatorname{supp}\mu\subset F,\quad
       \|\Gamma_\mu\|_\infty\leq1\bigr\}.
\end{equation}
Since the zero measure is admitted in the definition, the competing class for \eqref{capacity-definition} is
always nonempty. A compact set $F\neq\varnothing$ may nevertheless have zero capacity: for instance, $\operatorname{cap}(\{z_0\})=0$.  See Lemma~\ref{lem:box-capacity}.
Throughout the paper,
capacity is exclusively intended in the sense of \eqref{capacity-definition}.

\begin{definition}\label{D:reg}
Given a bounded open set $U\subset\R^{n+1}$, a point
$z_0=(x_0,t_0) \in\partial U$ is called \emph{$\cL$-regular} if, for every
$f\in C(\partial U)$, the Perron-Wiener-Brelot solution $H_f^U$ of
$\cL u=0$ in $U$ satisfies
\[
 \lim_{U\ni z\to z_0}H_f^U(z)=f(z_0).
\]
\end{definition}
The solution $H_f^U$ exists by \cite[Theorem~3.1 and Section~4.1]{Kogoj}.

\subsection{The Wiener-type criterion}

For $z_0\in\partial U$ and $\rho>0$, define the exterior kernel
superlevel set
\begin{equation}\label{exterior-superlevels}
E_\rho(z_0) =U^c\cap\left(\{\zeta=(\xi,\tau)\mid\tau<t_0,\quad
\Gamma(z_0;\zeta)\geq\rho^{-1}\}\cup\{z_0\}\right).
\end{equation}
For $a>1$, define the corresponding exterior kernel annuli by
\begin{equation}\label{wiener-annuli}
A_j(z_0) =U^c\cap\left(\{\zeta=(\xi,\tau)\mid\tau<t_0,\quad
a^j\leq\Gamma(z_0;\zeta)\leq a^{j+1}\}\cup\{z_0\}\right).
\end{equation}

For sufficiently small $\rho>0$ and sufficiently large $j$, these sets
are contained in an arbitrarily small neighborhood of $z_0$. Their
compactness follows after the pole is adjoined: away from the pole the
kernel is continuous, while the small-time estimates proved in
Lemmas~\ref{lem:local-geometry} and \ref{lem:box-capacity} rule out any
other accumulation point. Adding the pole does not change capacity, since as we have mentioned
$\operatorname{cap}(\{z_0\})=0$. 

We are ready to state the main result.

\begin{theorem}[Wiener's criterion]\label{thm:wiener}
Let $U\subset\R^{n+1}$ be a bounded open set, and let
$z_0\in\partial U$. The following conditions (i)-(iv) are equivalent:
\begin{itemize}
\item[(i)] $z_0$ is $\cL$-regular for $U$.
\item[(ii)] For every $a>1$ and some (equivalently, every) $j_0\in\mathbb N$,
\[
 \sum_{j=j_0}^{\infty}a^j\operatorname{cap}(A_j(z_0))=\infty.
\]
\item[(iii)] For every sufficiently small $\delta>0$,
\[
 \int_0^{\delta}\frac{\operatorname{cap}(E_\rho(z_0))}{\rho^2}\,d\rho
 =\infty.
\]
\item[(iv)] For every sufficiently small $\delta>0$,
\[
 \int_0^{\delta}
       \frac{\operatorname{cap}(E_{V(t)}(z_0))}{V(t)}\frac{dt}{t}
 =\infty.
\]
\end{itemize}
\end{theorem}

As in Wiener's original criterion \cite{Wie}, the conditions in Theorem~\ref{thm:wiener} are geometric and local in nature, in the sense that the divergence of the series in Theorem~\ref{thm:wiener}(ii) is not affected by deleting finitely many of its terms. Similarly, the divergence of the integrals in Theorem~\ref{thm:wiener}(iii) and (iv) is not affected by decreasing their positive upper limits. For instance, if for some nonempty bounded open set $O\subset\R^n$ and some $s_0>0$ the complement of $U$ contains the set 
\[
\big\{(e^{sB}x_0+\Delta_{\sqrt s}\eta,t_0-s)\mid\eta\in O,\ 0<s<s_0\big\},
\]
 then $\operatorname{cap}(E_\rho(z_0))\geq c\rho$ for all small $\rho>0$, and Theorem~\ref{thm:wiener}(iii) gives the regularity of $z_0$; this is the exterior cone criterion of Kogoj \cite[Theorem~6.2]{Kogoj}. The proof of Theorem~\ref{thm:wiener} follows the scheme of Evans and Gariepy \cite{EG} for the heat equation, i.e., $Q=I$, $B=0$ in \eqref{operator-definition}. As in \cite{EG} and \cite{GL88}, it rests on two ingredients: mean value formulas on the level sets of $\Gamma$, and a strong form of the Harnack inequality, which we now introduce. 

\subsection{The strong Harnack inequality} Theorem~\ref{thm:strong} below represents the main novelty of the present work. The sets $A_j(z_0)$ defined in \eqref{wiener-annuli} lie on level sets of $\Gamma$, and they contain regions exterior to every region with vertex at $z_0$ defined by the dilations \eqref{spacedil} in which the ordinary Harnack inequality controls a positive solution from below. Theorem~\ref{thm:strong} provides this control on a whole kernel ball, in terms of the average of the solution on an earlier time slice. 

\begin{figure}[htbp]
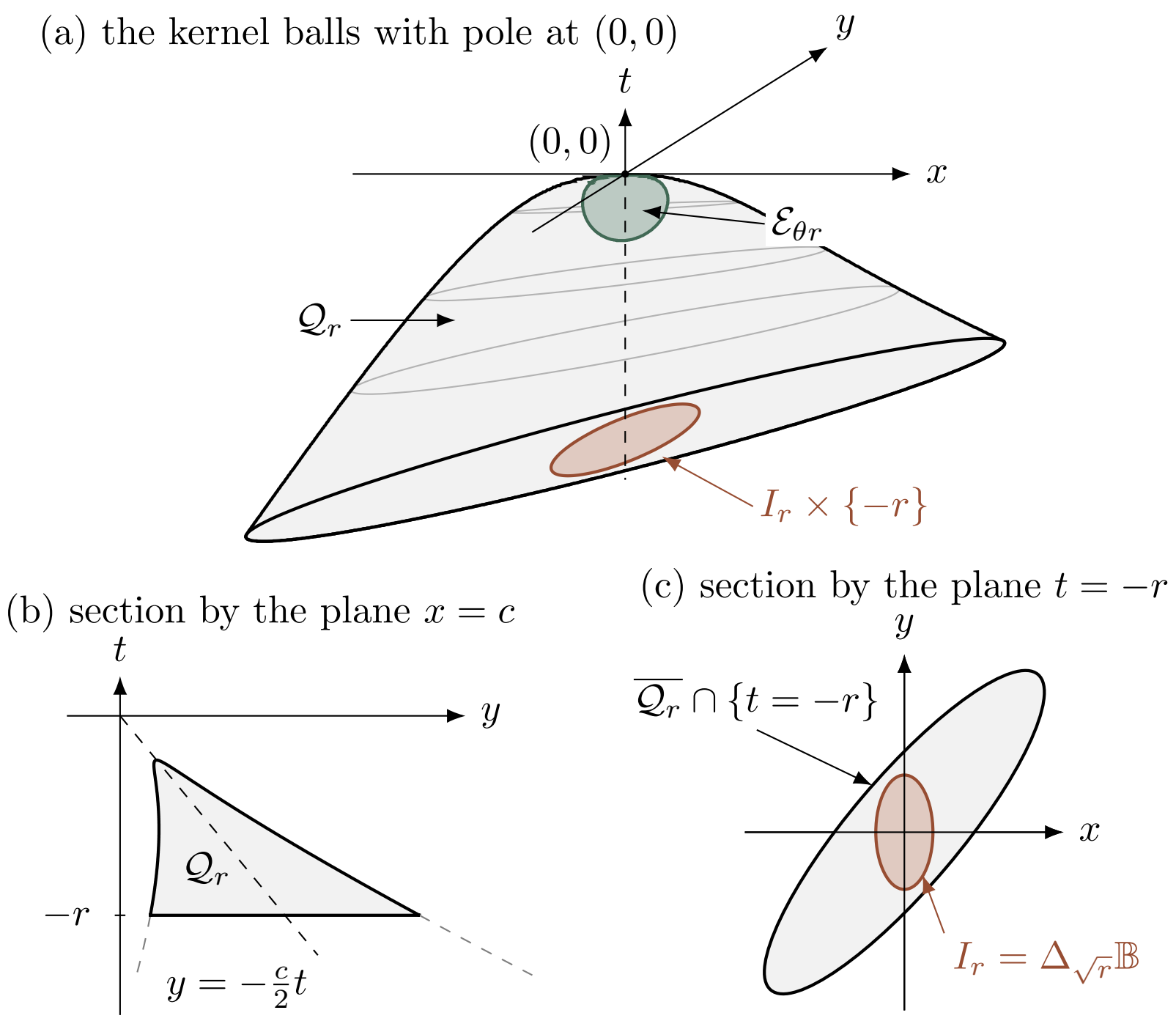

\centering

\caption{The geometry of Theorem~\ref{thm:strong} for the Kolmogorov equation \eqref{kolmo}, with pole at the origin. (a) An oblique view of the truncated kernel ball $\mathcal Q_r$, with its section at $t=-r$, of the smaller kernel ball $\mathcal E_{\theta r}$ on which the infimum in \eqref{strong-harnack} is taken, and of the set $I_r\times\{-r\}$ on which $u$ is averaged.  (b) The section of $\mathcal Q_r$ by the plane $x=c$. (c) The section of $\overline{\mathcal Q_r}$ at $t=-r$, and the ellipse $I_r$ inside it.\newline\mbox{}}
\label{fig:strong-harnack}
\end{figure}

We note that, while we state Theorem~\ref{thm:strong} at the origin, and accordingly denote by $\mathcal E_r$ the kernel balls of the operator \eqref{operator-definition} with pole at $(0,0)$, there is no loss of generality in doing so, since the case of an arbitrary pole $z_0$ follows by left translation with respect to the group law \eqref{Lie}. For $r>0$ and $\alpha>1$ we denote by
\[
\mathcal Q_r=\mathcal E_{\alpha r}\cap\{-r<t<0\}
\]
the truncated kernel ball and, with $\mathbb B\subset\Rn$ the Euclidean unit ball, by
\[
I_r=\Delta_{\sqrt r}\mathbb B
\]
the image of $\mathbb B$ under the anisotropic dilations \eqref{spacedil}, which are associated with the nilpotent principal part $\overline B$ of the matrix $B$ defined in \eqref{principal-B}. For the constants $\alpha$ and $\theta$ in Theorem~\ref{thm:strong}, and for $0<r\leq r_*$, the ellipsoid $I_r\times\{-r\}$ lies strictly inside the section of $\overline{\mathcal Q_r}$ at $t=-r$, and $\mathcal E_{\theta r}\subset\mathcal Q_r$. The geometry is illustrated in Figure \ref{fig:strong-harnack}.

\medskip

\begin{theorem}[Strong Harnack inequality]\label{thm:strong}
There are constants $\alpha>1$, $\theta\in(0,1)$, $C_H<\infty$,
and $r_*>0$, depending only on $Q,B$ in \eqref{canonical-QB}, with the following property.
For every $0<r\leq r_*$,
every $u\ge 0$ such that $\cL u=0$ in an open set containing $\overline{\mathcal Q_r}\setminus\{0\}$ satisfies
\begin{equation}\label{strong-harnack}
\displaystyle
 \frac1{|I_r|}\int_{I_r}u(x,-r)\dd x
       \leq C_H\inf_{\mathcal E_{\theta r}}u.
\end{equation}
The same constants apply at every pole by left translation with respect
to \eqref{Lie}. If $B=\overline B$, the estimate holds for every $r>0$.
\end{theorem}

\subsection{Description of the results} We next provide a brief description of our results. The proof of Theorem~\ref{thm:strong} follows the idea of Evans and Gariepy \cite[Lemma~3.2]{EG}, see also \cite[Theorem~1.4]{GL88}. In the region $|\Delta_{\sqrt{-t}}^{-1}x|\leq R$, where $R\geq1$ is a large constant depending only on $Q,B$, fixed in Theorem~\ref{thm:barrier}, and $\Delta_\lambda$ are the dilations in \eqref{spacedil}, a lower bound for $u$ in terms of its average on $I_r\times\{-r\}$ follows from the intrinsic Harnack inequality of \cite{LP}, combined near the pole with the barrier described below (Lemma~\ref{lem:core}). In the rest of the kernel ball one compares $u$ with $\Psi(\log((4\pi)^{n/2}V(r)\Gamma(0;\cdot)))$, where $\Psi$ is bounded and increasing, and $\Psi(0)=0$. This function vanishes on the boundary of $\mathcal E_r$, is bounded below on $\mathcal E_{\theta r}$, and is a subsolution wherever $|Q^{1/2}\nabla_x\log\Gamma|^2$ is dominated by a multiple of the derivative of $\log\Gamma$ along the vector field $\partial_t-\langle Bx,\nabla_x\rangle$. 

For the heat equation this happens \emph{exactly} outside a paraboloid with vertex at the pole, see \cite[(2.2)]{GL88}, and the maximum principle concludes the proof. For the operators \eqref{operator-definition} \emph{this is no longer true}. 

This obstruction can be seen already for \eqref{kolmo}. There this domination fails in a strip around the line $2sx-3y=0$, on which $m_s(0,(x,y))^2=x^2/(3s)$ is unbounded, and therefore the barrier has the wrong sign on a set which is not contained in any region of this form, see the blue shaded \emph{failure region} in Fig.~\ref{fig:heat-kolmogorov-comparison} below and the computation in Section~\ref{sec:model-geometry}. 

We remove this obstruction by subtracting from $\log\Gamma$ a bounded function $G(\Delta_{\sqrt{-t}}^{-1}x)$ whose derivative along the drift is positive where the diffusion gradient of $\Gamma$ vanishes, see Theorem~\ref{thm:barrier}. Since $G$ is bounded, the corrected barrier still vanishes outside $\mathcal E_r$ and is bounded below on a smaller kernel ball, and the rest of the argument goes through. In the rescaled variables the quadratic form of the kernel grows along the drift only at the rate of the diffusion gradient, and $G$ is a multiple of the logarithm of its ratio to a second quadratic form, which grows at a positive rate in every direction. The construction is explicit for \eqref{kolmo} in Section~\ref{sec:model-geometry}, and is carried out for a general pair $(Q,B)$ in Section~\ref{sec:strong-harnack}.

With Theorem \ref{thm:strong} in hand, we prove the sufficiency part of Theorem~\ref{thm:wiener}, i.e., that the divergence of the series in (ii) implies regularity; the necessity only uses the subadditivity of capacity and the bound \eqref{box-capacity} for the capacity of small boxes. As in \cite[Section~4]{GL88}, we take $z_0=0$ and let $V_F$ be the equilibrium potential of a compact set $F$ assembled from compact pieces of a sparse subfamily of the annuli $A_j(0)$. Let $A$ be one of these pieces, on which $\Gamma(0;\cdot)$ is comparable to $\rho^{-1}$, and let $r$ be comparable to the time radius of the kernel ball $\{\Gamma(0;\cdot)>\rho^{-1}\}$. If $\mu$ is the Riesz measure of $V_F$ and $\Omega$ is a kernel ball containing $\overline{\mathcal Q_r}\setminus\{0\}$, then $v=1-\Gamma(\mu|_{\Omega^c})$ is a nonnegative solution of $\cL v=0$ in $\Omega$, and $v\geq1-V_F$. 

Theorem~\ref{thm:strong} bounds $v$ from below on $A\subset\mathcal E_{\theta r}$ by $C_H^{-1}$ times the average of  $1-V_F$ on $I_r\times\{-r\}$. Integrating this bound against a measure on $A$ whose adjoint potential is bounded by $1$, and estimating $\mu(\Omega)$ by the mean value formulas of Section~\ref{sec:means}, we obtain in Lemma~\ref{lem:interaction} that $\rho^{-1}\operatorname{cap}(A)$, multiplied by this slice average, is bounded by a constant times the decrease of the normalized means of $V_F$ between two levels comparable to $\rho$. Along the sparse sequence of levels these decrements telescope, and their sum is bounded independently of the number of annuli. If the series in (ii) diverges, the slice averages of  $1-V_F$ must therefore tend to zero along a subsequence, and the semigroup identity together with \eqref{slice-mass} gives $V_F(0)\geq\delta$, with $\delta>0$ independent of the neighbourhood of $0$ in which $F$ is chosen. By \eqref{balayage-regularity}, this is the regularity of $0$.

In closing, we mention that in the companion paper \cite{GII} we use Theorem~\ref{thm:wiener} to prove a counterpart, for the class \eqref{operator-definition}, \eqref{kalman}, of the invariance result in \cite[Theorem~1.2]{GL88}: for the operators
\[
 \partial_t-\langle Bx,\nabla_x\rangle-\diver_x(\sigma a(x,t)\sigma^T\nabla_x),
 \]
 with $Q=\sigma\sigma^T$ and $a=[a_{ij}]$ smooth and uniformly positive definite, a boundary point is regular if and only if it is regular for the operator obtained by freezing $a$ at that point.

\subsection{Organization} The paper is organized as follows. In Section~\ref{sec:scaling} we collect the canonical representation of the pair $(Q,B)$, the small-time behavior of the Gramian and of the volume function, and the uniform intrinsic Harnack inequality, and in Section~\ref{sec:EG-prelim} we recall the Evans--Gariepy argument for the heat equation. In Section~\ref{sec:model-geometry} we examine the Kolmogorov model, and in Section~\ref{sec:strong-harnack} we prove Theorem~\ref{thm:strong} for a general pair $(Q,B)$. In Section~\ref{sec:means} we establish the potential-theoretic estimates, and in Section~\ref{sec:wiener} we prove Theorem~\ref{thm:wiener}.

\section{Preliminaries}
\label{sec:scaling}

In this section we collect some known properties of the operators \eqref{operator-definition} which will be used in the proofs of our main results. These are the canonical representation associated with the Kalman--H\"ormander condition, together with the anisotropic dilations determined by its nilpotent principal part; the resulting small-time behavior of the Gramian matrix $Q(t)$ and of the volume function $V(t)$; and the uniform intrinsic Harnack inequality, which will be used in Section~\ref{sec:strong-harnack} in the proof of the strong Harnack inequality.

A first basic observation that follows from a change of variable in the integral defining the Gramian \eqref{Qt} is 
\begin{equation}\label{mono}
Q(t+s) = Q(t) + e^{tB} Q(s) e^{tB^T},\ \ \ \ \ \ \ \ s, t >0.
\end{equation}
It immediately follows from \eqref{mono} and  Weyl's monotonicity theorem (see \cite[Cor.8.4.10,(iv)]{Bernstein})  that the  function $t\to V(t)$ defined in \eqref{intrinsic-distance-volume} is strictly increasing on $(0,\infty)$.

\subsection{Canonical representation and intrinsic dilations}

We recall that in the paper \cite{LP} it was shown that  \eqref{kalman} is equivalent to having $Q$ and $B$ in the following canonical form
\begin{equation}\label{canonical-QB}
 Q=\begin{pmatrix}
 Q_0&0\\[2pt]0&0
 \end{pmatrix},
 \qquad
 B=\begin{pmatrix}
 \star&\star&\star&\cdots&\star\\
 B_1&\star&\star&\cdots&\star\\
 0&B_2&\star&\cdots&\star\\
 \vdots&\ddots&\ddots&\ddots&\vdots\\
 0&\cdots&0&B_r&\star
 \end{pmatrix},
\end{equation}
where $Q_0 = Q_0^T$ is a strictly positive $p_0\times p_0$ matrix, and for $j=1,...,r$, $B_j$ is a $ p_j \times p_{j-1}$ matrix such that $\operatorname{Rank}(B_j) = p_j$,  with 
\[
p_0\ge p_1\ge ... \ge p_r \ge 1,\ \ \  \text{and}\ \ \ p_0 + p_1 + ... + p_r = n.
\]
The matrices designated by a $\star$ in the above expression of $B$ can be arbitrary. When $Q$ itself is invertible, and therefore $p_0 = n$,  we have $r=0$, and $B$ coincides with the $\star$ in the upper left corner.  We note that the canonical form depends on the pair $(Q,B)$: if $Q$ is changed, the canonical representation of $B$ changes accordingly. We refer to \cite[Sec. 7]{BG} for an extensive discussion of this aspect.
 
We emphasize that $B$ in \eqref{canonical-QB}
 need not be nilpotent and therefore, based on the results in \cite{Ku, Ku2, Kaw, LP}, there can be no global dilations attached to the PDE. 
However, it is possible to canonically assign a local scaling  to every pair $(Q,B)$.
The ranks $p_0, p_1,...,p_r$ of the matrices $Q_0, B_1,...,B_r$ in \eqref{canonical-QB} induce a corresponding stratification of the ambient space 
\[
\Rn=\R^{p_0}\oplus\R^{p_1}\oplus\cdots\oplus\R^{p_r}.
\]
We accordingly denote by 
\begin{equation}\label{deco}
x=\left(x^{(p_0)},x^{(p_1)},\ldots,x^{(p_r)}\right),\qquad x^{(p_j)}\in\R^{p_j},
\end{equation}
a generic point of $\Rn$, and define the non-isotropic dilations  
\begin{equation}\label{sdl}
\delta_\lambda(x,t) = \left(\lambda
x^{(p_0)},\lambda^3
x^{(p_1)},\ldots,\lambda^{2r+1}x^{(p_r)},\lambda^2t\right) = (e^{\log \la E} x,\la^2 t),
\end{equation}
where
\begin{equation}\label{canonical-dilations}
 E=\diag(I_{p_0},3I_{p_1},\ldots,(2r+1)I_{p_r}).
\end{equation}
We will denote by 
\begin{equation}\label{spacedil}
\Delta_\la x = e^{\log \la E} x,
\end{equation}
the space anisotropic dilations associated with the Kalman filtration.  
\begin{definition}\label{D:hd}
We call the number  
\begin{equation}\label{hd}
D = \tr E = p_0+3p_1+\ldots+(2r+1)p_r
\end{equation}
the \emph{local homogeneous dimension} of the pair $(Q,B)$. Throughout the paper we let\  $\nu=D/2$.
\end{definition}
The qualification ``local'' refers
to the fact that these dilations describe the geometry as time tends to zero,
even when the arbitrary $\star$ blocks in \eqref{canonical-QB} prevent the
operator $\cL$ in \eqref{operator-definition} from being globally homogeneous.
Note that such $\cL$ is homogeneous of degree two with respect to \eqref{sdl} if and only if $B = \overline B$, where $\overline B$ is the nilpotent matrix
obtained substituting all the $\star$ blocks in $B$ by zero matrices
\begin{equation}\label{principal-B}
 \overline B=\begin{pmatrix}
 0 & 0&0&\cdots&0\\
 B_1&0&0&\cdots&0\\
 0&B_2&0&\cdots&0\\
 \vdots&\ddots&\ddots&\ddots&\vdots\\
 0&\cdots&0&B_r&0
 \end{pmatrix}.
\end{equation}

For instance, the Kolmogorov model \eqref{mako} is 
in the form $B = \overline B$, and we have $p_0 = p_1 =1$. Consequently, according to \eqref{sdl} the dilations are
\[
\delta_\la(x,y,t) = (\la x,\la^3 y,\la^2 t),
\]  
and the homogeneous dimension is $D = 1+3=4$.
The following highly degenerate extension of \eqref{kolmo} was studied by Dym \cite{Dym}
\begin{equation}\label{badd}
\p_t u =  \p_{x_1x_1} u + x_1 \p_{x_2} u + . . . + x_{n-1} \p_{x_n} u,
\end{equation}  
which has dilations
\begin{equation}\label{badil}
\delta_\lambda(x,t) = \left(\lambda
x_1,\lambda^3
x_2,\ldots,\lambda^{2n-1}x_n,\lambda^2t\right).
\end{equation}
Its homogeneous dimension is $D=n^2$, and its Gramian is explicit, 
\[
Q(t)_{ij}=t^{i+j-1}/[(i-1)!(j-1)!(i+j-1)], \ \ \ \ \ 1\leq i,j\leq n,
\]
cf.\ \cite[Lemma~2.1]{Dym}.

Observe that \eqref{hd} implies that:
\begin{itemize}  
\item[(a)] $D\ge p_0+p_1+\ldots+p_r = n$, and $n\ge2$ when $Q$ is singular;
\item[(b)] when $Q$ is invertible, we have $r=0$, and therefore $D =n= p_0$;
\item[(c)] when $Q$ is not invertible, we have $D>n$. This immediately follows from \eqref{hd} since $p_0<n$.
\end{itemize}

Throughout the rest of this work, we assume that $Q$ and $B$ in \eqref{operator-definition} are in the canonical form \eqref{canonical-QB}. 

\subsection{Small-time behavior of the volume function}\label{S:smalltime}
We next recall a basic small time asymptotic of the volume function $V(t)$. Since we are going to use the Kalman filtration associated with $\overline B$, we first assume that $B = \overline B$ and observe that, in such case, $\delta_\lambda$ in \eqref{sdl} are group automorphisms with respect to the Lie group law \eqref{Lie}. This 
 follows from the commutation property established in \cite[p.288]{Ku2}, 
\begin{equation}\label{expdil}
e^{-\lambda^2 t \bar B}=\Delta_\lambda e^{-t \bar B}\Delta_{\lambda^{-1}},
\end{equation}
valid for any $\lambda>0$ and $t\in\R$.
We denote by 
\[
\overline V(t) = \sqrt{\det \overline Q(t)} = \left(\det \int_0^t e^{s\overline B} Q e^{s \overline B^T} ds\right)^{1/2},
\]
the volume function generated by $Q$ and $\overline B$. 
 
\begin{lemma}\label{kupcov}
Let $Q_*
 =
 \int_0^1
 e^{\sigma\overline B}
 Qe^{\sigma\overline B^T}\,d\sigma$.
For every $t>0$ we have
\begin{equation}\label{barQ}
\overline V(t)=t^{D/2}\sqrt{\det Q_*}.
\end{equation}
\end{lemma}

\begin{proof}
By \eqref{expdil}, $e^{t\sigma\overline B}=\Delta_{\sqrt t}e^{\sigma\overline B}\Delta_{\sqrt t}^{-1}$, and since by \eqref{canonical-QB}, $Q$ acts only on the first layer, we have $\Delta_{\sqrt t}^{-1}Q\Delta_{\sqrt t}^{-1}=t^{-1}Q$. The change of variables $s=t\sigma$ therefore gives $\overline Q(t)=\Delta_{\sqrt t}Q_*\Delta_{\sqrt t}$, and \eqref{barQ} follows from $\det\Delta_{\sqrt t}=t^{D/2}$.

\end{proof}

The explicit form of $Q_*$ was found in \cite{Ku2}. For instance, for the Kolmogorov operator \eqref{kolmo}, we have $r=1$, 
\[
Q_*=\begin{pmatrix}1&1/2\\1/2&1/3\end{pmatrix}
\]
and $\det\overline Q(t)=t^4/12$.

With the dilations $\Delta_\lambda$ in \eqref{spacedil}, consider the rescaled drift 
\[
B_s=s\Delta_{\sqrt{s}}^{-1}B\Delta_{\sqrt{s}}.
\]
 Its $(i,j)$ block is $s^{1+j-i}B_{ij}$, where $B_{ij}$ is the $(i,j)$ block of $B$ in \eqref{canonical-QB}, the block indices running from $0$ to $r$. In the canonical form $B_{ij}=0$ for $i>j+1$; for $i=j+1$ the factor is $s^0=1$, so the $(j+1,j)$ block of $B_s$ is the block $B_{j+1}$ in \eqref{canonical-QB} for every $s$; for $i\leq j$ the exponent $1+j-i$ is at least $1$. Hence
\begin{equation}\label{scaledB}
 B_s=\overline B+O(s)\quad(s\downarrow0).
\end{equation}
The asymptotic relation \eqref{scaledB} shows that $\overline B$ governs the small-time geometry.
Set 
\begin{equation}\label{hatQ}
\hat Q_s=\Delta_{\sqrt{s}}^{-1}Q(s)\Delta_{\sqrt{s}}^{-1}.
\end{equation}
 The change of variable $\sigma=sv$ in $Q(s)=\int_0^s e^{\sigma B}Qe^{\sigma B^T}\dd\sigma$ gives $Q(s)=s\int_0^1 e^{svB}Qe^{svB^T}\dd v$. Since $\Delta_{\sqrt{s}}^{-1}e^{svB}\Delta_{\sqrt{s}}=e^{vB_s}$ and $s\Delta_{\sqrt{s}}^{-1}Q\Delta_{\sqrt{s}}^{-1}=Q$, we obtain
\begin{equation}\label{scaledC}
 \hat Q_s=\int_0^1 e^{vB_s}Qe^{vB_s^T}\dd v.
\end{equation}
By \eqref{scaledB}, \eqref{scaledC} we infer 
\begin{equation}\label{hatQs}
\hat Q_s=Q_*+O(s).
\end{equation}
Since $Q_*$ is invertible, the same holds for the inverses.

We now justify the geometric assertions made after \eqref{kernel-ball-explicit}. By \eqref{kernel-ball-explicit} and \eqref{scaledC}, the spatial section of $\mathcal E_r(z_0)$ at time $t_0-s$, $0<s<r$, is the ellipsoid
\[
\frac14\langle\hat Q_s^{-1}\eta,\eta\rangle<\log(V(r)/V(s))
\]
in the variable 
\[
\eta=\Delta_{\sqrt s}^{-1}(e^{sB}x_0-y).
\]
Since $\lambda_{\max}(Q(s))\leq Cs$ and  $V(s)=s^\nu\sqrt{\det\hat Q_s}=s^\nu\sqrt{\det Q_*}\,(1+O(s))$, its Euclidean diameter, which is at most $4\sqrt{\lambda_{\max}(Q(s))\log(V(r)/V(s))}$, tends to zero as $s\to0^+$ and as $s\to r^-$. On the other hand, $\hat Q_s^{-1}$ stays bounded while $\log(V(r)/V(s))\to\infty$ as $s\to0^+$. Hence, for every $R>0$, the kernel ball contains points with $|\eta|>R$ arbitrarily close to its pole.

\subsection{The uniform intrinsic Harnack inequality}

We use the invariant Harnack inequality of Lanconelli and Polidoro
\cite[Theorems~5.1 and 5.2]{LP}.  As the authors explicitly observe, these results
extend and improve the Harnack inequality established by Garofalo and
Lanconelli in \cite[Theorem~1.1]{GL}, where the geometry determined by
the level sets of the fundamental solution was first introduced.

Given $R>0$ and $0<\sigma_2<\sigma_1<\sigma_0$, there is $R_0>R$, depending only on $Q,B,R,\sigma_0,\sigma_1,\sigma_2$, for which the following holds. Put
\[
 \mathcal C_r(R_0,\sigma_0)
 =\{(x,-s):0<s<\sigma_0r,
              |\Delta_{\sqrt r}^{-1}x|<R_0\}.
\]
There are constants $C<\infty$ and $r_0>0$, depending only on $Q,B,R,\sigma_0,\sigma_1$ and $\sigma_2$, such that every
nonnegative solution of $\cL u=0$ in
$\mathcal C_r(R_0,\sigma_0)$, $0<r\leq r_0$, satisfies
\begin{equation}\label{uniform-intrinsic-harnack}
 \sup_{|\Delta_{\sqrt r}^{-1}x|\leq R}u(x,-\sigma_1r)
 \leq C
 \inf_{|\Delta_{\sqrt r}^{-1}x|\leq R}u(x,-\sigma_2r).
\end{equation}
The restriction on $R_0$ cannot be dropped. Harnack inequalities of this sup/inf form hold for compact sets in the interior of the propagation set of the points where the infimum is taken, see \cite[Theorem~1.1]{CNP} for the case $B=\overline B$. For \eqref{kolmo}, in the rescaled variables, a backward admissible path that starts on the later slice and stays in $\mathcal C_r(R_0,\sigma_0)$ changes the second coordinate by at most a constant times $R_0(\sigma_1-\sigma_2)$ before it reaches time $-\sigma_1$; so when this quantity is small compared with $R$, part of the earlier slice is not reachable. With $R_0$ large, \eqref{uniform-intrinsic-harnack} follows from \cite[Theorem~5.2]{LP} by a finite chain of invariant cylinders along admissible paths, as in the proof of Lemma~\ref{lem:core}. The estimate at a general pole follows by left translation with respect to the group law \eqref{Lie}, which leaves $\cL$ invariant.
For the uniformity in $r$, note that under the change of variables
$x=\Delta_{\sqrt r}\xi$, $t=r\tau$, the drift matrix becomes
$B_r=r\Delta_{\sqrt r}^{-1}B\Delta_{\sqrt r}$.  By \eqref{scaledB},
$B_r\to\overline B$ as $r\downarrow0$. The rescaled operators are the operators $L_\lambda$, $\lambda=\sqrt r$, in the proof of \cite[Theorem~5.2]{LP}, where the constants are shown to be uniform for $0\leq\lambda\leq\lambda_0$; hence the constants in the
intrinsic Harnack inequality can be chosen uniformly for
$0<r\leq r_0$.

\subsection{The Evans--Gariepy strong Harnack inequality for the heat equation}
\label{sec:EG-prelim}

We close this preliminary section by recalling for the unfamiliar reader the argument of Evans and Gariepy for the strong Harnack inequality \cite[Lemma~3.2]{EG}, as presented in \cite[Section~2]{GL88}. We place the pole at the origin, write $s=-t>0$, and let
\[
 g(x,-s)=(4\pi s)^{-N/2}e^{-|x|^2/(4s)},\qquad
 L=\log\frac{g(x,-s)}{g(0,-r)}.
\]
Then $\{L>0\}$ is the heat ball whose lowest time is $-r$, and $L=0$ on its boundary, apart from the pole. Let $u\geq0$ be a solution of the heat equation, normalized so that its average on a fixed space ball in an earlier time slice is one. Since $|\nabla L|^2=|x|^2/(4s^2)$ and $-\partial_tL=|x|^2/(4s^2)-N/(2s)$, for every $\Theta>1$ one has
\[
 |\nabla L|^2\leq-\Theta\,\partial_tL
 \quad\text{if and only if}\quad
 |x|^2\geq\frac{2N\Theta}{\Theta-1}s,
\]
i.e., outside a paraboloid with vertex at the pole, see \cite[(2.2)]{GL88}. We fix $\Theta=16/9$, so that $2N\Theta/(\Theta-1)=32N/7$, and we denote by 
\[
P=\{(x,-s):s>0,\ |x|^2<(32N/7)s\}
\]
 the corresponding paraboloid. The strong Harnack inequality \cite[Theorem~1.4]{GL88} amounts to a lower bound $u\geq c_0>0$ on a smaller heat ball $\{L>\ell\}$, with $\ell>0$ fixed and $c_0$ independent of $r$. Its proof has three steps.

In the first step one bounds $u$ from below in $P$. Because of the time gap between the earlier slice and $P$, the ordinary Harnack inequality gives $u\geq c>0$ there, with $c$ independent of $r$; this is \cite[(2.33)]{GL88}.

In the second step one bounds $u$ from below in the part of the heat ball outside $P$, by comparing it with a function of $L$. In \cite{EG} and in \cite[(2.34)]{GL88} this function is $\arctan(L+16)-\arctan16$; any bounded increasing $\Psi$ with $\Psi(0)=0$ and $-\Psi''\leq C\Psi'$ would do. Here we use instead $v=\Psi(L)$, where $\Psi(\sigma)=8(1-e^{-\sigma/8})$, which simplifies the computation.
Since $-\Psi''=\Psi'/8$, we find
\[
 (\partial_t-\Delta)v
 =\Psi'(L)(\partial_t-\Delta)L-\Psi''(L)|\nabla L|^2
 =\Psi'(L)\left(\frac Ns-\frac7{32}\frac{|x|^2}{s^2}\right),
\]
and this is nonpositive outside $P$. Thus $v$ is a subsolution in the region which is not covered by the first step. Since $v=0$ on the heat sphere and $(c/8)v\leq c\leq u$ on the paraboloid, the maximum principle gives $u\geq(c/8)v$ outside the paraboloid, cf.\ \cite[(2.39)]{GL88}.

In the third step one observes that on the smaller heat ball $\{L>\ell\}$ the function $L$, and hence $v$, is bounded below by a positive constant, as in the passage from (2.39) to (2.40) in \cite{GL88}. Together with the first step, this gives the lower bound for $u$ on the whole smaller heat ball.

For the operators \eqref{operator-definition} the second step breaks down: the set where $v$ fails to be a subsolution need not be contained in a region $|\Delta_{\sqrt{-t}}^{-1}x|\leq R$, with $\Delta_\lambda$ as in \eqref{spacedil}, where the intrinsic Harnack inequality applies. In Section~\ref{sec:model-geometry} we exhibit this failure for the Kolmogorov model \eqref{kolmo}, and in Section~\ref{sec:strong-harnack} we correct the barrier for a general pair $(Q,B)$.

\section{The Kolmogorov model}
\label{sec:model-geometry}

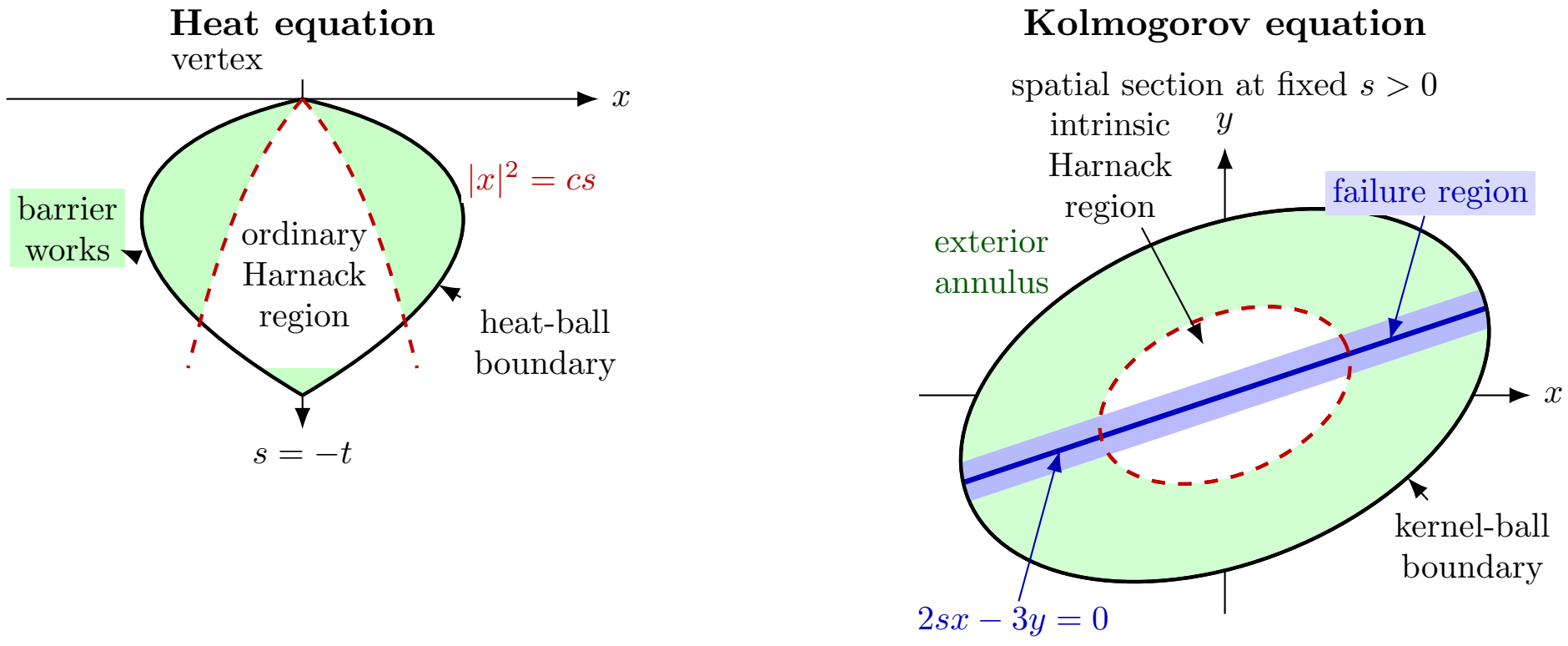
\begin{figure}[htbp]
\centering
\begin{tikzpicture}[scale=0.77,transform shape,font=\small,>=Latex]
  \begin{scope}[xshift=-4.8cm]
    \node[font=\bfseries] at (0,3.80) {Heat equation};
    \draw[->] (-3.05,3.05) -- (3.05,3.05) node[right] {$x$};
    \draw[->] (0,3.25) -- (0,-0.35) node[below] {$s=-t$};
    \node[anchor=south] at (-0.88,3.22) {vertex};

    \path[draw=black,thick,fill=gray!7]
      (0,3.05)
      .. controls (-1.95,2.60) and (-2.45,1.45) .. (0,0)
      .. controls (2.45,1.45) and (1.95,2.60) .. (0,3.05) -- cycle;

    \begin{scope}
      \clip (0,3.05)
        .. controls (-1.95,2.60) and (-2.45,1.45) .. (0,0)
        .. controls (2.45,1.45) and (1.95,2.60) .. (0,3.05) -- cycle;
      \fill[green!22] (-3.0,3.1) rectangle (3.0,-0.1);
      \fill[white]
        (0,3.05)
        .. controls (-0.62,2.35) and (-0.92,1.30) .. (-1.18,0.28)
        -- (1.18,0.28)
        .. controls (0.92,1.30) and (0.62,2.35) .. (0,3.05) -- cycle;
    \end{scope}

    \draw[black,thick]
      (0,3.05)
      .. controls (-1.95,2.60) and (-2.45,1.45) .. (0,0)
      .. controls (2.45,1.45) and (1.95,2.60) .. (0,3.05);
    \draw[red!75!black,thick,dashed]
      (0,3.05)
      .. controls (-0.62,2.35) and (-0.92,1.30) .. (-1.18,0.28);
    \draw[red!75!black,thick,dashed]
      (0,3.05)
      .. controls (0.62,2.35) and (0.92,1.30) .. (1.18,0.28);
    \path (0,0) .. controls (2.45,1.45) and (1.95,2.60) .. (0,3.05)
      coordinate[pos=0.28] (heatbdpt);

    \node[align=center,inner sep=1pt] at (0.02,1.22) {ordinary\\Harnack\\region};
    \node[align=center,fill=green!22,inner sep=2pt] (barworks) at (-2.42,1.72)
      {barrier\\works};
    \draw[->] (barworks) -- (-1.88,1.50);
    \node[align=center] (heatbound) at (2.50,0.52) {heat-ball\\boundary};
    \draw[->] (heatbound) -- (heatbdpt);
    \node[red!75!black,fill=white,inner sep=1.5pt] at (2.36,2.22) {$|x|^2=cs$};
  \end{scope}

  \begin{scope}[xshift=4.7cm]
    \node[font=\bfseries] at (0,3.80) {Kolmogorov equation};
    \node at (0,3.18) {spatial section at fixed $s>0$};
    \draw[->] (-3.15,0) -- (3.15,0) node[right] {$x$};
    \draw[->] (0,-2.25) -- (0,2.55) node[above] {$y$};

    \begin{scope}[rotate=22]
      \fill[green!18] (0,0) ellipse (2.85 and 1.72);
      \fill[white] (0,0) ellipse (1.35 and 0.82);
      \draw[black,thick] (0,0) ellipse (2.85 and 1.72);
      \draw[red!75!black,thick,dashed] (0,0) ellipse (1.35 and 0.82);
    \end{scope}

    \begin{scope}
      \clip[rotate=22] (0,0) ellipse (2.85 and 1.72);
      \fill[blue!30,opacity=0.90]
        (-3.3,-1.315) -- (3.3,0.885) -- (3.3,1.315) -- (-3.3,-0.885) -- cycle;
      \draw[blue!75!black,very thick]
        (-3.3,-1.1) -- (3.3,1.1);
    \end{scope}

    \begin{scope}[rotate=22]
      \draw[black,thick] (0,0) ellipse (2.85 and 1.72);
      \draw[red!75!black,thick,dashed] (0,0) ellipse (1.35 and 0.82);
    \end{scope}

    \node[align=center,inner sep=1pt] (intrh) at (-1.18,2.34)
      {intrinsic\\Harnack\\region};
    \draw[->] (intrh) -- (-0.22,0.52);
    \node[align=center,blue!70!black,fill=blue!15,inner sep=2pt]
      (fail) at (2.12,2.08) {failure region};
    \draw[->,blue!70!black] (fail) -- (1.70,0.567);
    \node[blue!70!black,anchor=north,inner sep=1.5pt] (bluelbl) at (-2.18,-2.12)
      {$2sx-3y=0$};
    \draw[->,blue!70!black] (bluelbl) -- (-1.70,-0.567);
    \node[align=center] (kerbound) at (2.55,-1.58) {kernel-ball\\boundary};
    \draw[->] (kerbound) -- (1.879,-0.847);
    \node[green!35!black,align=center] at (-2.40,1.38)
      {exterior\\annulus};
  \end{scope}
\end{tikzpicture}

\medskip
\begin{equation*}
\boxed{
\begin{gathered}
 b_r(X,-s)=\bigl[8\bigl(1-e^{-F_r(X,-s)/8}\bigr)\bigr]_+,\\
 F_r=\log\!\left(4\pi V(r)\Gamma_0(X,-s)\right)
              -G_s(X).
\end{gathered}}
\end{equation*}
\caption{Left: the heat barrier works outside the dashed curve,
while the ordinary Harnack inequality controls its interior.  Right: at a
fixed backward time, the blue strip in the Kolmogorov section crosses the
inner dashed ellipse and continues into the exterior region.  The
uncorrected heat barrier fails on this strip.  The boxed function is the
corrected Kolmogorov barrier constructed below.}
\label{fig:heat-kolmogorov-comparison}
\end{figure}

In this section we examine the Kolmogorov operator, for which the obstruction to the Evans--Gariepy barrier and the correction that removes it can be seen explicitly. Following
the convention of \cite{GL}, write
\[
 \mathcal K=\partial_{xx}+x\partial_y-\partial_t.
\]
With the sign convention of the present paper, the corresponding operator in the form \eqref{operator-definition} is $-\mathcal K=\partial_t-\partial_{xx}-x\partial_y$.
For the sake of brevity, we denote the spatial vector by
\[
 X=(x,y)^T.
\]
Figures \ref{fig:strong-harnack} and \ref{fig:heat-kolmogorov-comparison}
refer to this model, with pole at the origin. Its Gramian matrix and
backward kernel are
\begin{align}
 Q(s)&=\begin{pmatrix}s&s^2/2\\s^2/2&s^3/3\end{pmatrix},
                                           \label{figure-gramian}\\
 \Gamma_0(x,y,-s)
 &=\frac{\sqrt3}{2\pi s^2}
    \exp\left[-\frac{x^2}{4s}
                   -\frac{3(y-sx/2)^2}{s^3}\right].
                                           \label{figure-kernel}
\end{align}
In this example,
\[
 V(s)=\frac{s^2}{\sqrt{12}},\qquad
 m_s(0,X)^2=\frac{x^2}{4s}+\frac{3(y-sx/2)^2}{s^3}.
\]
We therefore define the backward ball of radius $r$ by
\begin{equation}\label{model-time-ball}
 \begin{aligned}
 \mathcal E_r
 &=\left\{(X,-s):s>0,\ \Gamma_0(X,-s)>
                         \frac{(4\pi)^{-1}}{V(r)}\right\}\\
 &=\left\{(X,-s):0<s<r,\ m_s(0,X)^2<2\log(r/s)\right\}.
 \end{aligned}
\end{equation}
The obstruction is visible on the line $2sx-3y=0$. With $E=\log\Gamma_0(x,y,-s)$, direct differentiation gives $(\partial_xE)^2=(2x/s-3y/s^2)^2$ and $(\partial_s+x\partial_y)E=(2x/s-3y/s^2)^2-2/s$, so that for $\Theta>1$ the adjusted Li--Yau inequality $(\partial_xE)^2\leq\Theta(\partial_s+x\partial_y)E$ fails in the strip $(2sx-3y)^2<\frac{2\Theta}{\Theta-1}s^3$. On the central line $m_s(0,(x,y))^2=x^2/(3s)$, and for small $s$ the strip contains points of the kernel ball \eqref{model-time-ball} which lie outside any fixed region $m_s(0,(x,y))^2<R^2$, where the intrinsic Harnack inequality gives a lower bound (Lemma~\ref{lem:core}). Figure~\ref{fig:heat-kolmogorov-comparison} compares this geometry with that of the heat equation.

\subsection{Construction of the correction}

We now construct the correction, working directly in the original variables.  Since
$t=-s$, the first-order vector field that replaces the heat time
derivative is
\[
 Y=-\partial_s-x\partial_y.
\]
The quadratic expression in the exponential of the kernel is
\[
 m_s(0,X)^2
 =\frac{x^2}{4s}+\frac{3(y-sx/2)^2}{s^3}.
\]
Direct differentiation gives
\[
 Y\bigl(m_s(0,X)^2\bigr)
 =\left(\frac{2x}{s}-\frac{3y}{s^2}\right)^2.
\]
This derivative vanishes on the line $2sx-3y=0$, exactly where the
uncorrected comparison function fails.  We therefore seek another
quadratic expression, with the same homogeneity as $m_s^2$, whose
$Y$-derivative is strictly positive away from the origin.  The most
general such expression is
\[
 H_s(x,y)=A_0\frac{x^2}{s}+A_1\frac{xy}{s^2}
                         +A_2\frac{y^2}{s^3}.
\]
One computes
\[
 sYH_s=(A_0-A_1)\frac{x^2}{s}
       +(2A_1-2A_2)\frac{xy}{s^2}
       +3A_2\frac{y^2}{s^3}.
\]
The simplest choice is to require
\[
 sYH_s=\frac{x^2}{s}+\frac{y^2}{s^3}.
\]
This determines
\[
 A_2=A_1=\frac13,
 \qquad A_0=\frac43,
\]
and hence
\begin{equation}\label{prototype-SH}
 H_s(x,y)=\frac{4x^2}{3s}+\frac{xy}{3s^2}
                         +\frac{y^2}{3s^3},
 \qquad
 sYH_s=\frac{x^2}{s}+\frac{y^2}{s^3}.
\end{equation}
The last expression is strictly positive whenever $(x,y)\ne(0,0)$.
In the variables $\eta=(x/\sqrt s,\,y/s^{3/2})$ one has $sY=\langle A\eta,\nabla_\eta\rangle$ with $A=\begin{pmatrix}1/2&0\\-1&3/2\end{pmatrix}$, and \eqref{prototype-SH} says that $H_s=\langle S\eta,\eta\rangle$, with $S=\begin{pmatrix}4/3&1/6\\1/6&1/3\end{pmatrix}$, and $sYH_s=|\eta|^2$.
The matrix construction in Section~\ref{sec:strong-harnack}, for a general pair $(Q,B)$, is the counterpart in arbitrary dimension of this elementary choice.

The function $\log H_s$ cannot be inserted by itself because it is
unbounded as $(x,y)$ approaches the origin.  The same singularity
occurs in $\log m_s^2$ and cancels in the difference. Set
\[
 \gamma_*=\inf_{s>0,\,X\ne0}\frac{H_s(X)}{m_s(0,X)^2},
 \qquad
 G_s(X)=c\log\frac{H_s(X)}{\gamma_*m_s(0,X)^2}.
\]
Here $c>0$ is a constant, which will be fixed as in Section~\ref{sec:strong-harnack}.
Direct comparison of the two quadratic forms gives constants
$c_1,c_2>0$, independent of $s$, such that
\[
 c_1m_s(0,X)^2\leq H_s(X)\leq c_2m_s(0,X)^2.
\]
Consequently $G_s\geq0$ and $G_s$ is bounded uniformly in $s$.
On the line $2sx-3y=0$, the derivative of $m_s^2$ under $Y$ vanishes,
whereas \eqref{prototype-SH} gives $YH_s>0$.  The derivative of $G_s$ thus supplies the term which is missing in the uncorrected comparison function, and this leads us to define the corrected logarithm
\[
 F_r(X,-s)=\log\!\left(4\pi V(r)\Gamma_0(X,-s)\right)-G_s(X),
\]
and the comparison function is
\[
 b_r(X,-s)=\left[8\left(1-e^{-F_r(X,-s)/8}\right)\right]_+.
\]
The differential inequality for this barrier will be verified in Section~\ref{sec:strong-harnack} for a general pair $(Q,B)$, which includes the Kolmogorov operator as a special case.

\section{The strong Harnack inequality}
\label{sec:strong-harnack}

In this section we prove Theorem~\ref{thm:strong}. We carry out the construction of Section~\ref{sec:model-geometry} for a general pair $(Q,B)$, using the canonical representation and the small-time estimates of Section~\ref{sec:scaling}. The matrix $A$ of the model is replaced by $A=\frac12E-\overline B$, and the three elementary equations which determined the coefficients of $H_s$ in \eqref{prototype-SH} are replaced by the matrix equation \eqref{general-Lyapunov-equation} below. With
$\eta=\Delta_{\sqrt{s}}^{-1}x$, differentiation gives
\[
 -s\bigl(\partial_s+\langle\overline Bx,\nabla_x\rangle\bigr)
 F(\Delta_{\sqrt{s}}^{-1}x)
 =\langle A\eta,\nabla_\eta F(\eta)\rangle,
 \qquad A=\frac12E-\overline B.
\]
Here $E$ is the diagonal matrix of the dilation weights defined in
\eqref{canonical-dilations}, and $\overline B$ is the nilpotent
principal part of $B$ defined in \eqref{principal-B}.

For later use, set
\begin{equation}\label{AandP}
 A=\tfrac12E-\overline B,\qquad P=\tfrac14 Q_*^{-1},\qquad
 \Phi(\eta)=\langle P\eta,\eta\rangle.
\end{equation}

The Gramian matrix $\overline Q(s)$ of the pair $(Q,\overline B)$, which appears in the definition of $\overline V$ before Lemma~\ref{kupcov}, satisfies $\overline Q(s)=\Delta_{\sqrt{s}}Q_*\Delta_{\sqrt{s}}$ and $\overline Q'(s)=Q+\overline B\,\overline Q(s)+\overline Q(s)\overline B^T$.
Differentiating at $s=1$ yields
\begin{equation}\label{lyapP}
 AQ_*+Q_*A^T=Q,\qquad
 A^TP+PA=\tfrac14Q_*^{-1}QQ_*^{-1}.
\end{equation}
Let 
\begin{equation}\label{W0}
W_0=\langle A\eta,\nabla_\eta\rangle
\end{equation}
 denote the derivative along the linear vector field $\eta\mapsto A\eta$. Then \eqref{lyapP} gives
\begin{equation}\label{Wphi}
 W_0\Phi=\frac14|Q^{1/2}Q_*^{-1}\eta|^2\geq0.
\end{equation}

\subsection{The kernel and the obstruction}

As in Section~\ref{sec:EG-prelim}, the proof has two parts: a lower bound in the region $|\Delta_{\sqrt{s}}^{-1}x|\leq R$, given by Lemma~\ref{lem:core}, and the corrected barrier of Theorem~\ref{thm:barrier} outside it. We begin by fixing the notation common to both.

As recalled after \eqref{kalman}, the Gramian \eqref{Qt} is positive definite.
Moreover, $\int_{\R^n}p(x,y,s)\,dy=1$ and
\begin{equation}\label{semigroup}
 \int_{\R^n}p(x,y,s)p(y,z,r)\,dy=p(x,z,s+r),
\end{equation}
see \cite[Lemma~2.4, (i) and (v)]{GTiumj}.
With $s=-t>0$, the backward kernel with pole at the origin is
\begin{equation}\label{kernel}
 \Gamma(0;(x,-s))=
 \frac{(4\pi)^{-n/2}}{V(s)}\exp\bigl(-m_s(0,x)^2\bigr).
\end{equation}
Here the first argument of $\Gamma$ is fixed at the origin, and $(x,-s)$ is the variable; as in \eqref{spacetime-kernel}, the pole of the fundamental solution is its second argument. A function $v$ is called a subsolution if $\cL v\leq0$.

\begin{theorem}[Corrected barrier]\label{thm:barrier}
There exist $s_0>0$, $R\geq1$, a bounded nonnegative function
$G\in C^\infty(\R^n\setminus\{0\})$, and $M=\sup G<\infty$, depending
only on $Q,B$, with the following properties. The function
$G$ is homogeneous of degree zero under ordinary
scalar multiplication. For $0<r<s_0$, set
\begin{equation}\label{barrier}
 \begin{split}
 F_r(x,-s)&=\log\!\left((4\pi)^{n/2}V(r)\Gamma(0;(x,-s))\right)
                       -G(\Delta_{\sqrt{s}}^{-1}x),\\
 b_r(x,-s)&=\left[8\left(1-e^{-F_r(x,-s)/8}\right)\right]_+.
 \end{split}
\end{equation}
Then the barrier $b_r$ is a distributional subsolution of $\cL$ on
\[
 \mathcal U=\{(x,-s):0<s<s_0,\ |\Delta_{\sqrt{s}}^{-1}x|>R\},
\]
and is a classical subsolution wherever it is positive. Moreover,
\begin{align}
 &0\leq b_r\leq8,\label{bounded}\\
 &b_r=0\quad\hbox{outside }\mathcal E_r,\label{outer}\\
 &b_r\geq 8(1-e^{-1/8})
       \quad\hbox{if }(4\pi)^{n/2}V(r)\Gamma(0;(x,-s))
       \geq e^{M+1},\label{inner}
\end{align}
where the last two assertions are restricted to $\mathcal U$.
If $B$ is its homogeneous principal part $\overline B$, one can take
$s_0=+\infty$.
\end{theorem}

Consider the matrix $\hat Q_s$ defined by \eqref{hatQ}. Since $Q(s)=\Delta_{\sqrt s}\hat Q_s\Delta_{\sqrt s}$, one has $m_s(0,x)^2=\frac14\langle\hat Q_s^{-1}\eta,\eta\rangle$, with $\eta=\Delta_{\sqrt s}^{-1}x$. By \eqref{hatQs} there are $0<\lambda\leq\Lambda$ such that $\lambda I\leq\hat Q_s\leq\Lambda I$ for $0<s<s_0$, hence
\[
 \frac{|\eta|^2}{4\Lambda}\leq m_s(0,x)^2\leq\frac{|\eta|^2}{4\lambda}.
\]
Therefore the region $|\Delta_{\sqrt s}^{-1}x|\leq R$ excluded in Theorem~\ref{thm:barrier}, where the lower bound is provided by Lemma~\ref{lem:core}, contains the set where $m_s(0,x)^2\leq R^2/(4\Lambda)$ and is contained in the set where $m_s(0,x)^2\leq R^2/(4\lambda)$. In this sense it is comparable with the regions $m_s(0,x)^2<R^2$ considered in the introduction and in Section~\ref{sec:model-geometry}. Since $0\leq G\leq M$, the correction changes a logarithmic kernel level by at most $M$.

\subsection{Some Gramian identities}

Let
\[
 \tilde Q(s)=\int_0^s e^{-vB}Qe^{-vB^T}\,dv.
\]
The identities from \cite[Section~5]{GT} used below are
\begin{align}
 Q(s)&=e^{sB}\tilde Q(s)e^{sB^T},\label{two-gramians}\\
 Q'(s)&=Q+BQ(s)+Q(s)B^T,\label{Qprime}\\
 \frac{d}{ds}\log\det Q(s)&=\tr(QQ(s)^{-1})+2\tr B=\tr(Q\tilde Q(s)^{-1}).\label{tracerelation}
\end{align}
In particular, direct differentiation of the backward kernel gives
\[
 \cL\log\Gamma(0;(x,-s))
 =\tr(QQ(s)^{-1})+\tr B
 -\frac14\langle Q(s)^{-1}QQ(s)^{-1}x,x\rangle.
\]

\subsection{The general form of the geometric obstruction}

We first examine the homogeneous principal operator, since the
obstruction is already present there. Direct differentiation gives
\begin{equation}\label{uncorrected-phase}
 \begin{aligned}
 s\cL\log\Gamma(0;(x,-s))
 &=\tr(QQ_*^{-1})
   -\frac14\left|Q^{1/2}Q_*^{-1}\eta\right|^2,\\
 s\left|Q^{1/2}\nabla_x\log\Gamma(0;(x,-s))\right|^2
 &=\frac14\left|Q^{1/2}Q_*^{-1}\eta\right|^2.
 \end{aligned}
\end{equation}
For the heat equation, the inequality $|\nabla_xE|^2\leq\Theta\,D_tE$, $E$ the logarithm of the fundamental solution, holds exactly outside a paraboloid with vertex at the pole. This was observed by Evans and Gariepy \cite[Lemma~3.2]{EG}, see \cite[(2.1)--(2.2)]{GL88}. For uniformly parabolic operators with variable coefficients, its intrinsic analogue is \cite[Lemma~2.1]{GL88}, whose proof rests on the small-time asymptotics of the fundamental solution established in \cite[Theorem~2.1]{GLma}, see also \cite[Lemma~B]{GLma}. For \eqref{operator-definition} and $\Theta>1$, the analogous inequality, which we call the adjusted Li--Yau inequality, is
\begin{equation}\label{general-adjusted-li-yau}
 \left|Q^{1/2}\nabla_x\log\Gamma(0;(x,-s))\right|^2
 \leq-\Theta\,\cL\log\Gamma(0;(x,-s)).
\end{equation}
As in Section~\ref{sec:EG-prelim}, this is the inequality that makes $\Psi(\log\Gamma)$ a subsolution.
By \eqref{uncorrected-phase}, it holds precisely where
\begin{equation}\label{general-adjusted-li-yau-region}
 \frac14\left|Q^{1/2}Q_*^{-1}\eta\right|^2
 \geq\frac{\Theta}{\Theta-1}
       \tr(QQ_*^{-1}).
\end{equation}
The exceptional region is therefore a tubular neighborhood of
$\ker(Q^{1/2}Q_*^{-1})$. The regions in which the
intrinsic Harnack estimate is applied, on the other hand, are the
sets $|\eta|<R$, equivalently sublevel sets of the function $\Phi$ in \eqref{AandP}. 
When $Q$ is singular, the tube in
\eqref{general-adjusted-li-yau-region} is unbounded along every nonzero
vector in $\ker(Q^{1/2}Q_*^{-1})$, while $\Phi$ tends to infinity along
the same ray. Thus, as for the Kolmogorov model in Section~\ref{sec:model-geometry}, the exceptional region of the adjusted Li--Yau inequality is not contained in any set $|\eta|<R$ with $R$ fixed.

The uncorrected Evans--Gariepy barrier for the kernel ball
$\mathcal E_r$ is
\[
 8\left[1-\left((4\pi)^{n/2}V(r)
 \Gamma(0;(x,-s))\right)^{-1/8}\right]
 \quad\text{on }\mathcal E_r.
\]
The chain rule gives
\begin{equation}\label{uncorrected-barrier-sign}
 \begin{aligned}
 &s\cL\left[8\left(1-\left((4\pi)^{n/2}V(r)
 \Gamma(0;(x,-s))\right)^{-1/8}\right)\right]\\
 &\quad=\left((4\pi)^{n/2}V(r)\Gamma(0;(x,-s))\right)^{-1/8}
 \left[\tr(QQ_*^{-1})
       -\frac7{32}\left|Q^{1/2}Q_*^{-1}\eta\right|^2\right].
 \end{aligned}
\end{equation}
Thus this expression has the desired sign when the diffusion gradient
is large, but it has the opposite sign wherever that gradient vanishes, that is, on the nonzero vectors of $\ker(Q^{1/2}Q_*^{-1})$. As observed above, this set meets every sufficiently large ellipsoid $\{\Phi=R^2\}$, and therefore the failure cannot be removed by taking the exterior radius $R$ larger. Such a difficulty does not arise for a uniformly parabolic operator. If in fact $Q>0$, then
\[
 \frac{|Q^{1/2}Q_*^{-1}\eta|^2}
      {\langle Q_*^{-1}\eta,\eta\rangle}
\]
has a positive minimum on the unit sphere, and the negative term in
\eqref{uncorrected-barrier-sign} eventually dominates in every
direction.

\subsection{Construction of the bounded correction}

To construct the correction, we note that the change of variables $v(x,t)=u(e^{-tB}x,t)$ transforms $\cL u=0$ into $\partial_tv-\tr(e^{tB}Qe^{tB^T}D_x^2v)=0$, see \cite[Lemma~2.1]{GLu}. In the variables of $v$ the derivative which replaces the heat time derivative is therefore $\partial_t$, and in the original variables it is the vector field
\begin{equation}\label{general-time-field}
 Y=\partial_t-\langle Bx,\nabla_x\rangle.
\end{equation}
Indeed, if $w(x,t)=F(e^{-tB}x,t)$, then
\[
 \partial_tw(x,t)=(YF)(e^{-tB}x,t).
\]
When $B=0$ we have $Y=\partial_t$, and therefore \eqref{general-time-field} is the counterpart, for the operator \eqref{operator-definition}, of the time derivative in the heat equation.

We next determine what the correction must accomplish. We look for it in the form $G_s(x)=G(\Delta_{\sqrt{s}}^{-1}x)$, where $G$ is the function in \eqref{G} below.  We note that the correction at $(x,-s)$ is a function of the original spatial variable and of the elapsed time, and that the anisotropic rescaling appears explicitly in its definition.

For the kernel ball $\mathcal E_r$, let
\[
 F_r(x,-s)=
 \log\!\left((4\pi)^{n/2}V(r)\Gamma(0;(x,-s))\right)-G_s(x)
\]
and set
\[
 \Psi(\sigma)=8(1-e^{-\sigma/8}).
\]
Where $F_r>0$, the chain rule gives
\begin{equation}\label{general-barrier-chain-rule}
 \cL\Psi(F_r)
 =e^{-F_r/8}\left\{
     \cL F_r+\frac18
     |Q^{1/2}\nabla_xF_r|^2\right\}.
\end{equation}
Our task is therefore to choose a bounded function $G_s$ for which the expression in braces is nonpositive; there is no need for $G_s$ to be large.  The
choice $G_s=0$ corresponds to the uncorrected barrier, and \eqref{uncorrected-barrier-sign} shows that it fails when
$Q^{1/2}Q_*^{-1}\eta$ is small.

We now return to the canonical variables and calculate which derivative
of $G_s$ enters \eqref{general-barrier-chain-rule}.  Recall that $t=-s$ and $\eta=\Delta_{\sqrt{s}}^{-1}x$. If $G_s(x)=G(\eta)$, differentiation at fixed $x$ gives
\[
 \partial_t\eta=\frac1{2s}E\eta.
\]
The drift term gives
\[
 -\langle Bx,\nabla_xG(\eta)\rangle
 =-\frac1s\langle B_s\eta,\nabla_\eta G(\eta)\rangle,
 \qquad
 B_s=s\Delta_{\sqrt{s}}^{-1}B\Delta_{\sqrt{s}}.
\]
Therefore
\begin{equation}\label{Y-in-canonical-variables}
 sY[G(\Delta_{\sqrt{s}}^{-1}x)]
 =\left\langle\left(\frac E2-B_s\right)\eta,
                         \nabla_\eta G(\eta)\right\rangle.
\end{equation}
For small $s$, $B_s=\overline B+O(s)$ by \eqref{scaledB}, and the first-order operator appearing in the limit is the operator $W_0$ in \eqref{W0}, with $A$ as in \eqref{AandP}. 

The logarithm of the limiting kernel is determined, apart from an additive constant, by the function $\Phi$ in \eqref{AandP}, and by \eqref{Wphi} its derivative $W_0\Phi$ is the quadratic expression coming from the diffusion gradient. It vanishes on the subspace
$\ker(Q^{1/2}Q_*^{-1})$.  Therefore the barrier cannot be repaired by using another function of $\Phi$, since every such function has zero $W_0$-derivative on the same subspace.

We instead look for a positive quadratic form
\[
 H(\eta)=\langle S\eta,\eta\rangle
\]
whose derivative along $W_0$ is positive in every nonzero direction.
Since
\[
 W_0H(\eta)=\langle(A^TS+SA)\eta,\eta\rangle,
\]
the simplest requirement is
\begin{equation}\label{general-Lyapunov-equation}
 A^TS+SA=I.
\end{equation}
The meaning of \eqref{general-Lyapunov-equation} is the following. Since the eigenvalues of $A$ are positive, see below, the flow $e^{vA}$ of $W_0$ carries every $\eta\neq0$ to infinity, and the solution of \eqref{general-Lyapunov-equation} is
\[
 H(\eta)=\int_0^\infty|e^{-vA}\eta|^2\dd v,
\]
the total squared length of the trajectory of $-W_0$ through $\eta$. Indeed, $H(e^{\tau A}\eta)=\int_{-\tau}^\infty|e^{-uA}\eta|^2\dd u$, and differentiation at $\tau=0$ gives $W_0H(\eta)=|\eta|^2$, which is \eqref{general-Lyapunov-equation}. Thus $H$ increases along the flow in every direction, at the rate $|\eta|^2$, whereas by \eqref{Wphi} the quadratic form $\Phi$ of the kernel increases only at the rate of the diffusion gradient. The function $G$ in \eqref{G} below is a multiple of the logarithm of $H/\Phi$. Being homogeneous of degree zero, it is bounded, and where the diffusion gradient is small its derivative along $W_0$ is close to $c|\eta|^2/H(\eta)$, which is bounded below by a positive constant, see \eqref{WG}. For the heat equation $A=\frac12I$, $H=4\Phi$ and $G$ is constant; for \eqref{kolmo}, $H$ is the quadratic form $H_s$ in \eqref{prototype-SH}. We next verify these facts and the properties needed in the barrier argument.

The matrix $\overline B$ is strictly lower block triangular, whereas
$E/2$ is block diagonal with diagonal blocks
\[
 \frac12I_{p_0},\ \frac32I_{p_1},\ldots,
 \frac{2r+1}{2}I_{p_r}.
\]
Consequently $A=E/2-\overline B$ is lower block triangular with these
same diagonal blocks. Its eigenvalues, with the appropriate
multiplicities, are
\[
 \frac12,\frac32,\ldots,\frac{2r+1}{2}.
\]
Hence $-A$ is asymptotically stable, and by the classical Lyapunov theorem \cite[Corollary~11.9.4 and Proposition~11.9.5]{Bernstein} the equation \eqref{general-Lyapunov-equation} has a unique solution, which is symmetric, positive definite, and given by
\begin{equation}\label{S}
 S=\int_0^\infty e^{-vA^T}e^{-vA}\dd v,
\end{equation}
that is, by the integral formula for $H$ above.

Let
\[
 \delta=\lambda_{\max}(S)^{-1}>0,\qquad
 \gamma_* =\min_{\eta\ne0}\frac{H(\eta)}{\Phi(\eta)}>0.
\]
For a constant $c>0$ to be fixed below, define
\begin{equation}\label{G}
 \boxed{\displaystyle
 G(\eta)=c\log\frac{H(\eta)}{\gamma_*\Phi(\eta)}.}
\end{equation}
Both quadratic forms are positive definite, so the quotient $H/\Phi$ is continuous and positive on the unit sphere, its minimum $\gamma_*$ is positive, and $G\geq0$. Since $H$ and $\Phi$ are homogeneous of degree two, $G$ is homogeneous of degree zero. It is smooth on $\R^n\setminus\{0\}$, where $H$ and $\Phi$ are smooth and positive, and it is bounded there, since it is bounded on the unit sphere.

Since $W_0H=|\eta|^2$ by \eqref{general-Lyapunov-equation}, and $W_0\Phi$ is given by \eqref{Wphi}, we obtain
\begin{equation}\label{WG}
 W_0G=c\left(\frac{|\eta|^2}{H(\eta)}
 -\frac{|Q^{1/2}Q_*^{-1}\eta|^2}
       {\langle Q_*^{-1}\eta,\eta\rangle}\right)
 \geq c\left(\delta-
 \frac{|Q^{1/2}Q_*^{-1}\eta|^2}
      {\langle Q_*^{-1}\eta,\eta\rangle}\right).
\end{equation}
The definition $\delta=\lambda_{\max}(S)^{-1}$ implies
\[
 H(\eta)\leq\lambda_{\max}(S)|\eta|^2,
 \qquad \frac{|\eta|^2}{H(\eta)}\geq\delta.
\]
In particular, $W_0G\geq c\delta/2$ whenever
\[
 \frac{|Q^{1/2}Q_*^{-1}\eta|^2}
      {\langle Q_*^{-1}\eta,\eta\rangle}\leq\frac\delta2.
\]
Hence the derivative of $G$ along $W_0$ is uniformly
positive on a neighbourhood of the set where the diffusion gradient
degenerates.

Homogeneity and smoothness on the unit sphere give two finite constants,
denoted by $C_D$ and $C_E$, such that
\begin{equation}\label{derivatives}
 |\tr(QD_\eta^2G)(\eta)|\leq\frac{C_D}{|\eta|^2},\qquad
 |Q^{1/2}\nabla_\eta G(\eta)|^2\leq\frac{C_E}{|\eta|^2}.
\end{equation}
These bounds follow directly by differentiating $G(r\eta)=G(\eta)$ and taking
the maxima of the resulting continuous functions on the unit sphere.

\subsection{The differential inequality}

\begin{proof}[Proof of Theorem~\ref{thm:barrier}]
We begin with some identities which hold for the operator \eqref{operator-definition}, without assuming $B=\overline B$. As before, we write $\eta=\Delta_{\sqrt{s}}^{-1}x$. Direct differentiation gives
\begin{equation}\label{exactkernel}
 \begin{aligned}
 s\cL\log\Gamma(0;(x,-s))
 &=\tr(Q\hat Q_s^{-1})+s\tr B
   -\frac14|Q^{1/2}\hat Q_s^{-1}\eta|^2,\\
 s|Q^{1/2}\nabla_x\log\Gamma(0;(x,-s))|^2
 &=\frac14|Q^{1/2}\hat Q_s^{-1}\eta|^2.
 \end{aligned}
\end{equation}
Indeed, $Q'(s)=Q+BQ(s)+Q(s)B^T$ implies
\[
 \begin{aligned}
 \partial_s\log\Gamma(0;(x,-s))
 &=-\tfrac12\tr\bigl(QQ(s)^{-1}\bigr)-\tr B\\
 &\quad+\tfrac14\langle Q(s)^{-1}QQ(s)^{-1}x,x\rangle
   +\tfrac12\langle Q(s)^{-1}Bx,x\rangle.
 \end{aligned}
\]
The remaining terms in $\cL\log\Gamma(0;(x,-s))$ are
\[
 \tfrac12\tr\bigl(QQ(s)^{-1}\bigr)
 +\tfrac12\langle Q(s)^{-1}Bx,x\rangle;
\]
the mixed terms cancel.
Combining these formulas gives \eqref{exactkernel}. We note the presence of the term $s\tr B$ in \eqref{exactkernel}.

For the correction, direct differentiation gives
\begin{equation}\label{exactG}
 s\cL[G(\Delta_{\sqrt{s}}^{-1}x)]
      =\left\langle\left(\frac E2-B_s\right)\eta,\nabla_\eta G(\eta)\right\rangle
       -\tr(QD_\eta^2G)(\eta).
\end{equation}
Indeed, differentiation with $x$ fixed gives
$\partial_t\eta=E\eta/(2s)$; the drift contributes $-B_s\eta/s$ and the
diffusion contributes $-\tr(QD_\eta^2G)/s$.

The scalar term
$\tr(Q\hat Q_s^{-1})+s\tr B$ extends continuously to
$s=0$. Choose $K\geq1$ such that its absolute value is at most $K$
for $0<s\leq1$, and now fix
\[
 c=\frac{2(K+3)}{\delta}
\]
in \eqref{G}. By \eqref{scaledB}--\eqref{scaledC}, uniformly in $\eta\ne0$,
\[
 \frac{|Q^{1/2}\hat Q_s^{-1}\eta|^2}
      {\langle Q_*^{-1}\eta,\eta\rangle}
 \longrightarrow
 \frac{|Q^{1/2}Q_*^{-1}\eta|^2}
      {\langle Q_*^{-1}\eta,\eta\rangle},
 \qquad
 \left\langle\left(\frac E2-B_s\right)\eta,\nabla_\eta G(\eta)\right\rangle
 \longrightarrow W_0G(\eta).
\]
The convergence is uniform in $\eta\neq0$, since the functions involved are homogeneous of degree zero and the unit sphere is compact. Choose $0<s_0\leq1$ so that for $0<s\leq s_0$,
\begin{align}
 \left|
 \frac{|Q^{1/2}\hat Q_s^{-1}\eta|^2
       -|Q^{1/2}Q_*^{-1}\eta|^2}
      {\langle Q_*^{-1}\eta,\eta\rangle}
 \right|&\leq\frac\delta4,\label{errors}\\
 \left|
 \left\langle\left(\frac E2-B_s\right)\eta,\nabla_\eta G(\eta)\right\rangle
 -W_0G(\eta)\right|&\leq1.\notag
\end{align}
If
\[
 \frac{|Q^{1/2}\hat Q_s^{-1}\eta|^2}
      {\langle Q_*^{-1}\eta,\eta\rangle}\leq\frac\delta4,
\]
then the corresponding quotient with $\hat Q_s$ replaced by
$Q_*$ is at most $\delta/2$, and hence
\begin{equation}\label{gooddrift}
 \left\langle\left(\frac E2-B_s\right)\eta,\nabla_\eta G(\eta)\right\rangle
 \geq K+2.
\end{equation}
Also fix $J\geq0$ such that
\[
 \left|\left\langle\left(\frac E2-B_s\right)\eta,
                  \nabla_\eta G(\eta)\right\rangle\right|\leq J
 \qquad(0<s\leq s_0).
\]
Let $p_*=\lambda_{\min}(P)>0$. Choose $R\geq1$ large enough that
\begin{equation}\label{Rchoice}
 \frac{C_D}{R^2}\leq1,\qquad \frac{C_E}{R^2}\leq1,
 \qquad \frac{\delta p_*R^2}{4}\geq2(K+J+2).
\end{equation}
We emphasize that none of these choices depends on $r$.

Combining \eqref{exactkernel} and \eqref{exactG}, for $|\eta|\geq R$,
\begin{equation}\label{LF}
 \begin{aligned}
 s\cL F_r
 &=\tr(Q\hat Q_s^{-1})+s\tr B
   -\frac14|Q^{1/2}\hat Q_s^{-1}\eta|^2\\
 &\qquad-
 \left\langle\left(\frac E2-B_s\right)\eta,\nabla_\eta G(\eta)\right\rangle
 +\tr(QD_\eta^2G)(\eta).
 \end{aligned}
\end{equation}
If the quotient in the preceding paragraph is at most $\delta/4$,
equations \eqref{gooddrift} and \eqref{Rchoice} give
\[
 s\cL F_r
 \leq-1-\frac14|Q^{1/2}\hat Q_s^{-1}\eta|^2.
\]
If that quotient is larger than $\delta/4$, then
\[
 \frac14|Q^{1/2}\hat Q_s^{-1}\eta|^2
 \geq\frac{\delta p_*|\eta|^2}{4}\geq2(K+J+2),
\]
and therefore
\[
 s\cL F_r
 \leq K+J+1-\frac14|Q^{1/2}\hat Q_s^{-1}\eta|^2
 \leq-1-\frac18|Q^{1/2}\hat Q_s^{-1}\eta|^2.
\]
Thus throughout the exterior region,
\begin{equation}\label{strict}
 s\cL F_r
 \leq-1-\frac18|Q^{1/2}\hat Q_s^{-1}\eta|^2.
\end{equation}
Furthermore, scaling the spatial gradients and using \eqref{derivatives},
\begin{equation}\label{energy}
 s\left|Q^{1/2}\nabla_x
 \left[\log\!\left((4\pi)^{n/2}V(r)\Gamma(0;(x,-s))\right)-G(\eta)\right]\right|^2
 \leq\frac12|Q^{1/2}\hat Q_s^{-1}\eta|^2+2.
\end{equation}

The chain rule, with the sign convention of $\cL$, now gives
\begin{align*}
 &s\cL\left[8\left(1-e^{-F_r/8}\right)\right]\\
 &\quad=e^{-F_r/8}\\
 &\qquad\times\left\{s\cL F_r
 +\frac{s}{8}\left|Q^{1/2}\nabla_x
 F_r\right|^2\right\}\\
 &\quad\leq-e^{-F_r/8}\\
 &\qquad\times\left(\frac34+\frac1{16}
 |Q^{1/2}\hat Q_s^{-1}\eta|^2\right)\\
 &\quad<0.
\end{align*}
It remains to show that the positive part $b_r$ is a subsolution in the sense of distributions. To see this, we apply the chain rule to the increasing convex functions $j_\epsilon(\lambda)=(\lambda+\sqrt{\lambda^2+\epsilon^2})/2$ composed with
$8(1-e^{-F_r/8})$, and pass to the locally uniform limit as
$\epsilon\downarrow0$. The additional diffusion term is
$-j_\epsilon''$ times the square of the diffusion gradient and is
nonpositive; hence it preserves the subsolution inequality.

The inequalities $0\leq G\leq M$ show that $F_r\leq0$ outside
$\mathcal E_r$ and that $F_r\geq1$ wherever
$(4\pi)^{n/2}V(r)\Gamma(0;(x,-s))\geq e^{M+1}$. They imply
\eqref{bounded}--\eqref{inner}. If $B=\overline B$, then $B_s=\overline B$,
$\hat Q_s=Q_*$, and $\tr B=0$ for every $s>0$. All the above estimates then hold for every $s>0$.
\end{proof}

\subsection{Completion of the proof}

\begin{corollary}[Exterior propagation]\label{cor:comparison}
Use the constants of Theorem~\ref{thm:barrier}, and let $0<r<s_0$. Put
\[
 \mathcal D_r=\mathcal E_r\cap
 \{(x,-s):|\Delta_{\sqrt{s}}^{-1}x|>R\}.
\]
Let $u$ be a nonnegative classical solution of $\cL u=0$ in a
neighbourhood of $\mathcal D_r$, continuous on its boundary away from
the pole. If, for some $m_0>0$,
\[
 u\geq m_0\quad\hbox{on }\{|\Delta_{\sqrt{s}}^{-1}x|=R,
                        \ (x,-s)\in\mathcal E_r\},
\]
then
\[
 u\geq m_0(1-e^{-1/8})
\]
where $(4\pi)^{n/2}V(r)\Gamma(0;(x,-s))\geq e^{M+1}$ and
$|\Delta_{\sqrt{s}}^{-1}x|>R$.
\end{corollary}

\begin{proof}
We compare $u$ with $(m_0/8)b_r$. On the part of the boundary of $\mathcal D_r$ contained in $\partial\mathcal E_r$ the function $b_r$ vanishes, and on the part where $|\Delta_{\sqrt{s}}^{-1}x|=R$ we have $(m_0/8)b_r\leq m_0$. For $0<\epsilon<r$ the sections of $\mathcal D_r$ in the strip $-r\leq t\leq-\epsilon$ are bounded, and the parabolic comparison principle gives $u\geq(m_0/8)b_r$ in $\mathcal D_r\cap\{t<-\epsilon\}$. Letting $\epsilon\downarrow0$, and using \eqref{inner}, we obtain the conclusion.
\end{proof}

For a general pole $z_0=(x_0,t_0)$ we replace $s=-t$ and $x$ by $s=t_0-t$ and $\xi=x-e^{sB}x_0$. In the variables $(\xi,t-t_0)$ the operator $\cL$ keeps its form, since the time derivative of $e^{sB}x_0$ cancels the corresponding drift term, and moreover $p(x_0,x,s)=p(0,\xi,s)$. Therefore all the constants are independent of the pole.

Corollary~\ref{cor:comparison} is the analogue of
\cite[Lemma 3.2, Step 4]{EG}. Its hypothesis on
$|\Delta_{\sqrt{s}}^{-1}x|=R$ will be obtained
from an average on an earlier time slice, by means of the intrinsic Harnack inequality recalled in Section~\ref{sec:scaling}.

\paragraph{Uniform rescaling.}

Recall that $\nu=D/2$, with $D$ as in \eqref{hd}.
Fix $r_*>0$ with $r_*\leq\min(1,s_0)$, where $s_0$ is supplied by
Theorem~\ref{thm:barrier}. For $0<\rho\leq r_*$ use the matrices
$B_\rho$ from \eqref{scaledB}; at $\rho=0$ use their limit $\overline B$.
Set
\[
 \mathcal Q_\rho(v)=\int_0^v e^{wB_\rho}Qe^{wB_\rho^T}\dd w.
\]
The rescaling identities are
\begin{equation}\label{double-scaling}
 v\Delta_{\sqrt v}^{-1}B_\rho \Delta_{\sqrt v}=B_{\rho v},\qquad
 \mathcal Q_\rho(v)=\Delta_{\sqrt v}\hat Q_{\rho v}\Delta_{\sqrt v},
\end{equation}
with $\hat Q_0=Q_*$. The backward kernel of the rescaled operator is
\[
 \Gamma_\rho(0;(z,-v))
 =\frac{(4\pi)^{-n/2}}{\sqrt{\det\mathcal Q_\rho(v)}}
       \exp\left(-\frac14z^T\mathcal Q_\rho(v)^{-1}z\right).
\]
For $\rho>0$ it satisfies
\begin{equation}\label{kernel-scaling}
 \Gamma_\rho(0;(z,-v))=\rho^\nu
                  \Gamma(0;(\Delta_{\sqrt\rho}z,-\rho v)).
\end{equation}
On $[0,r_*]$, the matrix $\hat Q_\sigma$ is continuous and positive
definite. Hence there are constants $c_-,c_+,K_0>0$ such that, uniformly for
$0\leq\rho\leq r_*$ and $0<v\leq1$,
\begin{equation}\label{gaussian-envelope}
 c_-v^{-\nu}\exp\left(-K_0|\Delta_{\sqrt v}^{-1}z|^2\right)
 \leq \Gamma_\rho(0;(z,-v))\leq c_+v^{-\nu}.
\end{equation}

\paragraph{The lower bound in the region controlled by the intrinsic Harnack inequality.}
\label{sec:core}

Let $\mathbb B=\{y\in\R^n:|y|<1\}$. For the $\gamma$ in
the following lemma, put
\[
 \mathcal Q_{\rho,\gamma}
 =\{(z,-w):0<w<1,\ \Gamma_\rho(0;(z,-w))>\gamma\}.
\]

\begin{lemma}[Estimate from a slice average]\label{lem:core}
For every $R\geq1$ there are $\gamma\in(0,1)$ and $c_0>0$,
depending only on $Q,B,R,r_*$, such that every nonnegative classical
solution $u$ of $\cL_\rho u=0$ in an open set containing $\overline{\mathcal Q_{\rho,\gamma}}\setminus\{0\}$ satisfies
\begin{equation}\label{core-average}
 u(\Delta_{\sqrt s}\eta,-s)\geq
    \frac{c_0}{|\mathbb B|}\int_{\mathbb B}u(y,-1)\dd y,
 \qquad |\eta|\leq R,\quad 0<s\leq\frac12.
\end{equation}
The constant is independent of $\rho$ and $s$.
\end{lemma}

\begin{proof}
Increase $R$, if necessary, so that the exterior barrier estimates
hold for $|\eta|\geq R$. Extend the function $G$ used there to a smooth
function $\widetilde G$ on $\R^n$, equal to $G$ for $|\eta|\geq R$
and satisfying $0\leq\widetilde G\leq M$. Such an extension is obtained
by multiplying $G$ by a smooth radial cutoff that vanishes near zero.
For a small constant $\varepsilon>0$, put
\[
 F=\log\frac{\Gamma_\rho(0;(x,-s))}{\varepsilon}
             -\widetilde G(\Delta_{\sqrt s}^{-1}x),
 \qquad b=8(1-e^{-F/8})_+.
\]
The estimates \eqref{strict}--\eqref{energy}, after normalization,
give $\cL_\rho b\leq0$ for $|\Delta_{\sqrt s}^{-1}x|\geq R$.
After decreasing $r_*$ if necessary, their constants are uniform for $0\leq\rho\leq r_*$. In the remaining region, smoothness of
$\widetilde G$ and the Gramian bounds give
\[
 |s\cL_\rho F|+s|Q^{1/2}\nabla_xF|^2\leq C,
 \qquad e^{-F/8}\leq C\varepsilon^{1/8}s^{\nu/8}
 \quad(F>0).
\]
Consequently, with $\beta=\nu/8>0$,
\begin{equation}\label{core-barrier-error}
 \cL_\rho b\leq C\varepsilon^{1/8}s^{\beta-1},
 \qquad 0<s<\tfrac12.
\end{equation}
As in the proof of Theorem~\ref{thm:barrier}, this inequality is preserved, in the sense of distributions, by taking the positive part. Therefore
\begin{equation}\label{core-time-correction}
 w=b-\frac C\beta\varepsilon^{1/8}
                      (2^{-\beta}-s^\beta)
\end{equation}
satisfies $\cL_\rho w\leq0$ on this time interval.

Set $b_*=8(1-e^{-1/8})$. On
$|\Delta_{\sqrt s}^{-1}x|\leq R$, $0<s\leq1/2$, the lower bound in
\eqref{gaussian-envelope} is bounded below by a positive constant
independent of $\rho$. Choose $\varepsilon$ small enough that
$F\geq1$ there and
$C\varepsilon^{1/8}2^{-\beta}/\beta\leq b_*/2$. Then
\begin{equation}\label{core-corrected-lower}
 w\geq b_*/2
 \quad\text{on }|\Delta_{\sqrt s}^{-1}x|\leq R,\quad0<s\leq1/2.
\end{equation}

It remains to compare at the earlier time $-1/2$. The sections
\[
 H_\rho=\{(x,-1/2):\Gamma_\rho(0;(x,-1/2))\geq\varepsilon\}
\]
are contained in a fixed compact set. The Kalman condition and
the uniform Gramian bounds connect $\overline{\mathbb B}\times\{-1\}$
to these sections in a fixed positive time, with paths in a fixed
bounded region. A finite chain of intrinsic Harnack inequalities
\cite[Theorem~5.2]{LP} therefore gives
\[
 \inf_{H_\rho}u\geq c_1\frac1{|\mathbb B|}
                              \int_{\mathbb B}u(y,-1)\,dy.
\]
The constant does not depend on $\rho$, since the normalized Gramians are uniformly positive definite and the matrices $B_\rho$ belong to a compact family. Choose $\gamma<\varepsilon$ smaller than the
uniform minimum of $\Gamma_\rho$ on the compact region containing
these chains. The chains then lie in $\mathcal Q_{\rho,\gamma}$.
Since $u$ solves the equation in a neighbourhood of $\overline{\mathbb B}\times\{-1\}$, the first cylinders of these chains may be placed below the slice $\{t=-1\}$. 

On the lateral boundary of
$\{\Gamma_\rho>\varepsilon,\ 0<s<1/2\}$ one has $w\leq0$,
and on its bottom section $w\leq8$. The comparison principle,
first up to time $-\delta$ and then with $\delta\downarrow0$, gives
\[
 u\geq\frac{c_1w}{8|\mathbb B|}
                         \int_{\mathbb B}u(y,-1)\,dy.
\]
Combining this with \eqref{core-corrected-lower} proves the lemma.
\end{proof}

\paragraph{Conclusion of the proof.}\label{sec:strong}

Recall that the averaging ellipsoid in Theorem~\ref{thm:strong} is
\begin{equation}\label{level-and-slice}
 I_r=\Delta_{\sqrt r}\mathbb B.
\end{equation}
Since $V$ is strictly increasing, the elapsed times $s=t_0-t$ in the kernel ball \eqref{kernel-ball-explicit} fill the interval $(0,r)$. When $Q=I_n$ and $B=0$, \eqref{kernel-ball-explicit} is the heat ball $0<s<r$, $|x_0-y|^2<2ns\log(r/s)$.

For every $r>0$ the closure of $\mathcal E_r$, with the pole adjoined, is compact. Indeed, its elapsed times lie in $[0,r]$, by the small-time estimates for the Gramian its spatial sections shrink to the pole as $s\downarrow0$, and at $s=r$ the section reduces to a single point. When $B\neq\overline B$ the kernel balls are not dilates of one another, and this is why Theorem~\ref{thm:strong} is stated for $0<r\leq r_*$.

\begin{proof}[Proof of Theorem \ref{thm:strong}]
Take $R,M,s_0$ from Theorem~\ref{thm:barrier} and use this $R$ in
Lemma~\ref{lem:core}. Choose
$\alpha>1$ so large that $c_+\alpha^{-\nu}<\gamma$.
Here $c_+$ is chosen so that the upper bound in \eqref{gaussian-envelope} holds whenever $\rho v\leq s_1$, for a fixed $s_1>0$; reduce $r_*$ so that $\alpha r_*\leq s_1$ and $r_*/2<s_0$. This does not affect $\gamma$, since Lemma~\ref{lem:core} remains valid for the smaller range of $\rho$. For $0<r\leq r_*$ define
$u_r(z,\tau)=u(\Delta_{\sqrt r}z,r\tau)$. Then $\cL_r u_r=0$.
By \eqref{kernel-scaling}, the rescaled outer level is
\[
 r^\nu\frac{(4\pi)^{-n/2}}{V(\alpha r)}
 \leq c_+\alpha^{-\nu}<\gamma.
\]
Consequently $u_r$ is a solution in an open set containing $\overline{\mathcal Q_{r,\gamma}}\setminus\{0\}$.
We can therefore apply Lemma~\ref{lem:core} there. Since the change of variables leaves the average over the bottom slice unchanged, we denote this average by
\[
 A_r=\frac1{|I_r|}\int_{I_r}u(x,-r)\dd x.
\]
Lemma~\ref{lem:core}, rescaled to the original variables, gives
\begin{equation}\label{physical-core}
 u(x,-s)\geq c_0 A_r
 \quad\text{when }0<s\leq r/2,\quad
                   |\Delta_{\sqrt s}^{-1}x|\leq R.
\end{equation}

We now apply the corrected barrier to the kernel ball
$\mathcal E_{r/2}$.  This ball is contained in
$\{-r/2<t<0\}$ and in $\mathcal E_{\alpha r}$.  Thus the equation
holds in a neighbourhood of its closure away from the pole, and
\eqref{physical-core} gives the lower bound $m_0=c_0A_r$ on
the portion $|\Delta_{\sqrt s}^{-1}x|=R$ of its boundary. If $A_r=0$,
\eqref{strong-harnack} is immediate.
Otherwise, Corollary~\ref{cor:comparison} with $r/2$ in place of $r$ gives
\[
 u\geq c_0(1-e^{-1/8})A_r
\]
where $(4\pi)^{n/2}V(r/2)\Gamma(0;(x,-s))\geq e^{M+1}$ and
$|\Delta_{\sqrt s}^{-1}x|>R$.
On the complementary part, where $|\Delta_{\sqrt s}^{-1}x|\leq R$, the
same lower bound follows from \eqref{physical-core}.
The small-time bounds for $V$ allow us to choose $0<\theta<1/2$,
independently of $r$, so that
\[
 \frac{V(r/2)}{V(\theta r)}>e^{M+1},
 \qquad 0<r\leq r_*.
\]
Indeed, $V(s)$ is bounded above and below by fixed positive multiples
of $s^\nu$ for small $s$.  Put
\[
 C_H=\bigl[c_0(1-e^{-1/8})\bigr]^{-1}.
\]
If $(x,-s)\in\mathcal E_{\theta r}$, then
\[
 (4\pi)^{n/2}V(r/2)\Gamma(0;(x,-s))
 >\frac{V(r/2)}{V(\theta r)}>e^{M+1}.
\]
Thus the preceding two lower bounds cover all of
$\mathcal E_{\theta r}$. This proves its inclusion in
$\mathcal Q_r$ and \eqref{strong-harnack}.
For the bottom-section inclusion, note that in the rescaled variables $\Gamma_r>\gamma$ on $\overline{\mathbb B}\times\{-1\}$, by the choice of $\gamma$ in the proof of Lemma~\ref{lem:core}, whereas the rescaled outer level is smaller than $\gamma$; hence $I_r\times\{-r\}$ lies in the bottom section of $\mathcal Q_r$.

The case of an arbitrary pole follows from the translation described after Corollary~\ref{cor:comparison}. Finally, when $B=\overline B$ we have $B_r=\overline B$ for every $r>0$, so that the rescaled operator does not depend on $r$, the intrinsic Harnack inequality holds at every scale, and the barrier has no restriction on time. The same proof then applies for every $r>0$.
\end{proof}

\section{Potential-theoretic estimates}\label{sec:means}

In this section and in Section~\ref{sec:wiener} we use the strong Harnack inequality in the Wiener argument. We first establish the normalized mean identities, extend them to bounded potentials, and derive the one-annulus mass estimate; we then relate this estimate to capacity, and in Section~\ref{sec:wiener} we carry out the summation over the exterior annuli.

\subsection{Normalized means and bounded potentials}

In this subsection we establish the mean value formulas of \cite{GL88} for the operators \eqref{operator-definition}. Since $\tr B$ need not vanish, these formulas contain an additional term, which we keep track of by means of the Gramian identities \eqref{two-gramians}--\eqref{tracerelation}. Throughout this section the parameter $\rho$ is the reciprocal of the height of the kernel. 

In this section only, abbreviate $\Gamma(0;(x,-s))$ by
$\Gamma(x,-s)$, and write
\[
 \Omega_\rho=\{(x,-s):s>0,\ \Gamma(x,-s)>\rho^{-1}\}.
\]
The kernel balls are among these sets, since
\begin{equation}\label{radius-level-conversion}
 \mathcal E_r=\Omega_{(4\pi)^{n/2}V(r)}.
\end{equation}
Choose $\rho_*>0$ so that these sets are bounded for
$0<\rho\leq\rho_*$. Their maximum ages
\[
 T_\rho=\sup\{s:(x,-s)\in\Omega_\rho\}
\]
satisfy $T_{(4\pi)^{n/2}V(r)}=r$ and tend to zero with $\rho$.
Indeed, $\Gamma(x,-s)\leq\Gamma(0,-s)$, the latter
function is strictly decreasing, and it tends to infinity as $s\downarrow0$.
The Gramian asymptotics in Section~\ref{sec:scaling} also show that
the spatial sections shrink to the pole as $s\downarrow0$.
Define
\begin{align}
 K_\rho
 &=|Q^{1/2}\nabla_x\log\Gamma|^2
   +(\tr B)\log(\rho\Gamma),\label{adjusted-kernel}\\
 D_\rho&=\log(\rho\Gamma),\qquad
 \mathcal M_\rho u
 =\frac1\rho\int_{\Omega_\rho}K_\rho u\,dx\,ds .
       \label{adjusted-mean}
\end{align}
Both kernels are extended by zero off $\Omega_\rho$.

\begin{lemma}[The inhomogeneous solid identity]\label{lem:solid}
Suppose $u$ is $C^2$ on a neighborhood of
$\overline{\Omega_{\rho_*}}$, including the pole. Then
\begin{equation}\label{solid-identity}
 u(0)=\mathcal M_\rho u+
          \int J_\rho(x,-s)\,\cL u(x,-s)\,dx\,ds,
 \qquad 0<\rho\leq\rho_*,
\end{equation}
where $J_\rho\geq0$ is given explicitly below. In particular,
\begin{align}
 \mathcal M_\rho1&=1, \label{mean-normalization}\\
 \frac{d}{d\rho}\mathcal M_\rho u
   &=-\frac1{\rho^2}
       \int D_\rho\,\cL u\,dx\,ds. \label{mean-derivative}
\end{align}
\begin{equation}\label{mean-total-variation}
 \frac1\rho\int |K_\rho|\,dx\,ds
 \leq 1+2|\tr B|T_\rho .
\end{equation}
\end{lemma}

\begin{proof}
The Gramian identities \eqref{two-gramians}--\eqref{tracerelation}
give
\[
 \cL^*\Gamma=0,
 \qquad
 \cL^*(\Gamma-\lambda^{-1})=-(\tr B)\lambda^{-1}
\]
away from the pole. For a.e. $\lambda$ the level set $\{\Gamma=\lambda^{-1}\}$ is smooth away from the pole, by Sard's theorem. As in the mean value formulas on level sets of \cite{GL} (see also \cite[Theorem~1.5]{GL88}), Green's identity on $\{\Gamma>\lambda^{-1}\}\cap\{t<-\epsilon\}$, followed by $\epsilon\downarrow0$, gives the usual surface mean formula together
with the additional term
$(\tr B)\lambda^{-1}\int_{\{\Gamma>\lambda^{-1}\}}u$.
The top face tends to $u(0)$ because the Gaussian has unit spatial
mass and its Gramian tends to zero.

Average this identity over $0<\lambda<\rho$ and use the coarea formula.
This gives \eqref{adjusted-kernel}--\eqref{solid-identity}, with
\begin{equation}\label{solid-remainder}
 J_\rho(x,-s)=
 \Gamma(x,-s)-\frac{1+\log(\rho\Gamma(x,-s))}{\rho}
\end{equation}
on $\Omega_\rho$, and $J_\rho=0$ outside it. Since
$a-1-\log a\geq0$ for
$a\geq1$,
\begin{equation}\label{remainder-domination}
 0\leq J_\rho\leq\Gamma,
 \qquad 0\leq D_\rho\leq\rho\Gamma.
\end{equation}
Substitution of $u=1$ proves \eqref{mean-normalization}, while
$\partial_\rho J_\rho=\rho^{-2}D_\rho$ proves
\eqref{mean-derivative}.

The only possibly negative part of $K_\rho$ is the term containing
$\tr B$, and
\[
 \rho^{-1}\int D_\rho\,dx\,ds\leq T_\rho.
\]
Together with \eqref{mean-normalization}, this proves
\eqref{mean-total-variation}.
\end{proof}

Formula \eqref{mean-derivative} is the counterpart of \cite[Theorem~1.6]{GL88}. Since $D_\rho\geq0$, the sign of $\frac{d}{d\rho}\mathcal M_\rho u$ is that of $-\cL u$. When $\tr B<0$ the kernel $K_\rho$ may change sign, and $\mathcal M_\rho$ is not a positive average. The mean value property therefore takes the following form: if $\cL u\geq0$ in $\Omega_{\rho_*}$, then $\rho\mapsto\mathcal M_\rho u$ is nonincreasing on $(0,\rho_*]$ and tends to $u(0)$ as $\rho\downarrow0$. The limit follows from \eqref{solid-identity}, since $0\leq J_\rho\leq\Gamma$ and $\int_{\Omega_\rho}\Gamma\,dx\,ds\leq T_\rho$.  Identities \eqref{solid-identity} and \eqref{potential-solid} are the counterparts, for the operators \eqref{operator-definition}, of the mean value formulas of Fabes and Garofalo \cite{FG} and Garofalo and Lanconelli \cite[Theorem~1.5]{GL88} for parabolic equations with variable coefficients; for operators of Kolmogorov type see \cite{GL}. These formulas were extended by Lanconelli and Pascucci \cite[Theorems~1.5 and~1.6]{LP99} to a class of hypoelliptic operators; their kernel contains the same term with the divergence of the drift, here $\tr B$.

\subsubsection{Extension to bounded potentials}

For a finite positive measure $\mu$ of compact support, write
$V_\mu=\Gamma\mu$ and suppose $0\leq V_\mu\leq1$. The standard
time regularization of a potential gives smooth functions
$0\leq u_\epsilon\uparrow V_\mu$ such that
$\cL u_\epsilon=f_\epsilon\geq0$ and
$f_\epsilon\,dz\to\mu$ vaguely. One may take
\[
 u_\epsilon(x,t)=\int\Gamma(x,t;\xi,\tau)
 F_\epsilon(t-\tau)\,d\mu(\xi,\tau),
\]
where $F_\epsilon$ is a smooth increasing approximation of the
Heaviside function supported away from the origin. The semigroup
identity shows that $u_\epsilon=\Gamma(f_\epsilon\,dz)$.

Applying Lemma~\ref{lem:solid} to $u_\epsilon$ and passing to the
limit gives
\begin{align}
 V_\mu(0)
 &=\mathcal M_\rho V_\mu+
       \int J_\rho\,d\mu,\label{potential-solid}\\
 \frac{d}{d\rho}\mathcal M_\rho V_\mu
 &=-\frac1{\rho^2}\int D_\rho\,d\mu .
                                      \label{potential-derivative}
\end{align}
Indeed, $0\leq J_\rho\leq\Gamma(0;\cdot)$ and
$0\leq D_\rho\leq\rho\Gamma(0;\cdot)$, while the total masses
$u_\epsilon(0)$ converge to $V_\mu(0)$. These bounds allow us to pass to the limit in both terms. Consequently the means are nonincreasing,
converge to $V_\mu(0)$ as $\rho\downarrow0$, and satisfy
\begin{equation}\label{potential-mean-bounds}
 -|\tr B|T_\rho\leq\mathcal M_\rho V_\mu
                         \leq1+|\tr B|T_\rho .
\end{equation}

This is the only place where we need to pass from smooth functions to equilibrium potentials.

\subsubsection{Quantitative mass and one-annulus capacity estimates}

\begin{corollary}[{Quantitative decrement; cf.\ \cite[Corollary~1.2]{GL88}}]\label{cor:mass-decrement}
For $0<2\rho\leq\rho_*$,
\begin{equation}\label{mass-decrement}
 \mathcal M_\rho V_\mu-\mathcal M_{2\rho}V_\mu
 \geq\frac{1-\log2}{2\rho}
     \mu\bigl(\{\Gamma\geq1/\rho,\ t<0\}\bigr).
\end{equation}
\end{corollary}

\begin{proof}
Integrate \eqref{potential-derivative} from $\rho$ to $2\rho$.
For $\rho<\lambda\leq2\rho$ and $\Gamma\geq1/\rho$,
$\log(\lambda\Gamma)\geq\log(\lambda/\rho)$. Restricting the positive
integrand to this set and using
\[
 \int_\rho^{2\rho}\frac{\log(\lambda/\rho)}{\lambda^2}\,d\lambda
 =\frac{1-\log2}{2\rho}
\]
gives \eqref{mass-decrement}. Since $\log(\lambda/\rho)>0$ for $\lambda>\rho$, the mass of $\mu$ carried by $\{\Gamma=1/\rho\}$ is included in this estimate.
\end{proof}

For a compact set $A$, write its equilibrium potential as
$V_A=V_{\mu_A}$, where $0\leq V_A\leq1$ and
$\operatorname{supp}\mu_A\subset A$.

\begin{remark}[On the need for a separate interaction estimate]
For one fixed bounded potential, the decrements telescope. For example,
with $\rho_j=2^{-j}\rho_*$,
\[
 \frac{1-\log2}{2}\sum_{j=1}^\infty\frac1{\rho_j}
           \mu(\{\Gamma\geq1/\rho_j,\ t<0\})
 \leq V_\mu(0)-\mathcal M_{\rho_*}V_\mu
 \leq1+|\tr B|T_{\rho_*}.
\]
This proves the measure estimate for bounded potentials, including the passage from smooth functions to their Riesz measures. Before the estimate can be summed over annuli, however, one has to take into account that the restriction of the equilibrium measure of a large compact set to one annulus is generally not the equilibrium measure of that annulus, so
their masses cannot be identified. Lemma~\ref{lem:interaction}
resolves this difficulty by comparing the annular capacity with a
mean decrement of one fixed equilibrium potential. The comparison also involves the average of $1-V$ on an earlier time slice, and it is to this average that the strong Harnack inequality will be applied.
\end{remark}

\subsection{Small-time geometry and capacity}\label{sec:local-geometry}

As before, $\nu=D/2$ with $D$ as in \eqref{hd}. Write
$\Gamma_0(z)=\Gamma(0;z)$, recall the sets $\Omega_\rho$ of Section~\ref{sec:means}, and set
\[
 \kappa=\frac{1-\log2}{2}>0.
\]
All level parameters will be small enough that the sets in question
lie in a fixed bounded time interval. We continue to denote by $V(s)=\sqrt{\det Q(s)}$ the volume function, and by $V_F$ the equilibrium potential of a compact set $F$.

\begin{lemma}[Local volume and slice estimates]\label{lem:local-geometry}
There are positive constants $s_*,b_-,b_+,\delta,C$ such that:
\begin{align}
 \begin{split}
 V(s)&=\sqrt{\det Q_*}\,s^\nu(1+O(s)),\\
 \frac{V'(s)}{V(s)}&=\frac\nu s+O(1);
 \end{split}\label{local-volume}\\
 b_-r^{-\nu}&\leq\Gamma_0(x,-r)\leq b_+r^{-\nu}
       \quad(x\in\overline I_r,\ 0<r\leq s_*);\label{slice-levels}\\
 \int_{I_r}p(0,x,r)\,dx&\geq\delta,\qquad
 |I_r|\sup_{I_r}p(0,\cdot,r)\leq C.\label{slice-mass}
\end{align}
In particular, $V(2s)\leq Cs^\nu\leq C'V(s)$ for $2s\leq s_*$. 
The closed sets
\[
 \{\Gamma_0\geq\rho^{-1},\ t<0\}\cup\{0\}
\]
are compact for sufficiently small $\rho$, and shrink to $0$ as
$\rho\downarrow0$.
\end{lemma}

\begin{proof}
Taking determinants in \eqref{scaledC} gives
$V(s)=s^\nu\sqrt{\det \hat Q_s}$. The entries of $B_s$ are polynomials in $s$,
so the integral defining $\hat Q_s$ is analytic at $s=0$.
It follows that $\hat Q_s=Q_*+O(s)$ with bounded derivative, and hence both
assertions in \eqref{local-volume} hold. The first is the local volume
asymptotic discussed in \cite[Appendix]{BG}; see also
\cite[(3.14) and Remark 3.1]{LP} and the underlying dilation
structure in \cite{Ku,Ku2}. 

For $x=\Delta_{\sqrt r}\eta$, the kernel is
\[
 \Gamma_0(x,-r)=r^{-\nu}
       \frac{(4\pi)^{-n/2}}{\sqrt{\det \hat Q_r}}
       \exp\left(-\frac14\langle\hat Q_r^{-1}\eta,\eta\rangle\right).
\]
Uniform upper and lower bounds for $\hat Q_r$ and $\hat Q_r^{-1}$ on a small
closed interval prove \eqref{slice-levels}. Since
$|I_r|=|\mathbb B|r^\nu$, integration gives
\eqref{slice-mass}, reducing $\delta$ if necessary so that
$0<\delta<1$.

To check the last assertion, $\Gamma_0(x,-s)\leq\Gamma_0(0,-s)$
bounds the elapsed time $s$ on $\Omega_\rho$ by a quantity which tends to zero with $\rho$. On a level set,
the Gramian bounds then give
\[
 |\Delta_{\sqrt s}^{-1}x|^2
       \leq C\log(C\rho s^{-\nu}).
\]
Since the right-hand side must be nonnegative, this forces $s\leq C'\rho^{1/\nu}$.
Consequently every layer satisfies
$|x^{(j)}|\leq C\rho^{(2j+1)/(2\nu)}$, by maximizing
$s^{2j+1}\log(C\rho s^{-\nu})$ over that interval.
These estimates show that the superlevel sets are bounded and shrink
to the origin. Their closedness after the pole is adjoined follows
from the continuity of the kernel away from the pole.
\end{proof}

\subsubsection{Capacity, balayage, and the dual test measure}
\label{sec:capacity-framework}

We use the kernel capacity defined in \eqref{capacity-definition}. When $\tr B\neq0$ the capacity of a compact set is not simply the mass of its equilibrium measure, and for this reason we state explicitly the potential-theoretic facts and the normalization which will be used below. The capacity is monotone and countably subadditive on compact sets. To see the latter, given a countable covering of a compact set by compact sets, one restricts an admissible measure to a subordinate disjoint measurable partition, and observes that each restriction is still admissible.

Although the condition in \eqref{capacity-definition} is a global bound, it suffices to verify it in a finite time strip containing $F$. Indeed, choose
$S<\min_F t\leq\max_F t<T$. Before $S$ the potential vanishes.
If it is bounded by $1$ on $S\leq t\leq T$, then for $t>T$ the
semigroup identity gives
\[
 \Gamma\mu(x,t)=\int_{\R^n}p(x,y,t-T)\Gamma\mu(y,T)\,dy\leq1,
\]
since $\int p(x,y,t-T)\,dy=1$. Thus the two admissibility conditions are equivalent, and no estimate of the volume function for large times is needed.

We shall use the following standard consequences of balayage theory,
in the form employed in \cite{GL88}. For a compact set $F$, its
equilibrium potential has the representation
$V_F=\Gamma\mu_F$, where $\mu_F$ is supported on $F$,
$0\leq V_F\leq1$, and $V_F=1$ in $\operatorname{int}F$.
The potentials are monotone and countably subadditive, and they
dominate every admissible potential supported on $F$. Moreover, if
\[
 C_R=\overline I_R\times[-R,0],\qquad F_R=C_R\setminus U,
\]
then
\begin{equation}\label{balayage-regularity}
 0\text{ is regular for }U
 \quad\Longleftrightarrow\quad
 \lim_{R\downarrow0}V_{F_R}(0)>0.
\end{equation}
Finally, for decreasing compact sets $F_j\downarrow F$ in a fixed
bounded region,
\begin{equation}\label{compact-continuity}
 V_{F_j}(z)\downarrow V_F(z),\qquad z\notin F.
\end{equation}
The harmonic-space structure for the operators \eqref{operator-definition} is given in \cite[Theorem~3.1]{Kogoj}, see also \cite{CC}. In particular, by \cite[Section~4.2]{Kogoj} we have $V_F=1$ in $\operatorname{int}F$, and \eqref{balayage-regularity} is \cite[Proposition~4.1]{Kogoj} with $B_\lambda=\overline I_\lambda$. The remaining facts follow as in the uniformly parabolic setting. Only the displayed statements are used below.

\subsubsection{The trace weight and duality}

Fix $T_*>0$ so that all sets and means used subsequently lie in
the slab $|t|\leq T_*$. Set
\[
 m(t)=e^{t\tr B},\qquad
 m_-=e^{-|\tr B|T_*},\quad m_+=e^{|\tr B|T_*},\quad
 \chi=m_+/m_-.
\]
We stress that no weight appears in the definition \eqref{capacity-definition} of capacity.

\begin{lemma}[Weighted mass and a local dual measure]
\label{lem:weighted-capacity}
For every compact $A\subset\{|t|\leq T_*\}$,
\begin{equation}\label{capacity-mass-comparison}
 \chi^{-1}\operatorname{cap}(A)\leq\mu_A(A)
                                      \leq\operatorname{cap}(A).
\end{equation}
There is a positive measure $\omega_A$ supported on $A$ such that
\begin{equation}\label{dual-test}
 \omega_A(A)\geq\chi^{-1}\operatorname{cap}(A),\qquad
 \Phi_A(z):=\int_A\Gamma(\zeta;z)\,d\omega_A(\zeta)\leq1
               \quad(|t|\leq T_*).
\end{equation}
\end{lemma}

\begin{proof}
We first note that
\begin{equation}\label{integrated-mass}
 \int_{\R^n}p(x,y,a)\,dx=e^{-a\tr B},
\end{equation}
see \cite[Lemma~2.4, (iii)]{GTiumj}.
Thus, for $T>\max\{\tau:(\xi,\tau)\in A\}$,
\[
 \int_{\R^n}m(T)V_A(x,T)\,dx
       =\int_A m(\tau)\,d\mu_A(\xi,\tau)=:c_m(A).
\]
If $\mu$ is admissible in \eqref{capacity-definition}, domination
and Tonelli's theorem show that
$\int_A m\,d\mu\leq c_m(A)$. Since $\mu_A$ itself is admissible,
it maximizes this weighted mass. Therefore
\begin{equation}\label{weighted-extremum}
 m_-\operatorname{cap}(A)\leq c_m(A)
                  \leq m_+\mu_A(A),\qquad
 \mu_A(A)\leq\operatorname{cap}(A).
\end{equation}
These inequalities give \eqref{capacity-mass-comparison}.

To construct the dual measure, consider the reversed evolution
\[
 \widehat{\cL}=-\partial_t-\tr(QD_x^2)+\langle Bx,\nabla_x\rangle.
\]
It is the operator of the class \eqref{operator-definition} with drift matrix $-B$, written in the reversed time variable. The pair $(Q,-B)$ still satisfies \eqref{kalman}. Write $p_{-B}$ for
the transition kernel with that drift, $p_B=p$, and
\[
 \widehat\Gamma((x,t);(y,s))=
 \begin{cases}p_{-B}(x,y,s-t),&s>t,\\0,&s\leq t.\end{cases}
\]
The Gramian identity
\eqref{two-gramians} gives the exact equality
\begin{equation}\label{weighted-duality}
 m(t)p_B(x,y,t-s)=m(s)p_{-B}(y,x,t-s),\qquad s<t.
\end{equation}
Let $\widehat V_A=\widehat\Gamma\widehat\mu_A$ be the equilibrium
potential for this reversed evolution, and put
$\widehat c_m(A)=\int_A m\,d\widehat\mu_A$.
We claim
\begin{equation}\label{weighted-capacity-duality}
 c_m(A)=\widehat c_m(A).
\end{equation}

Let $F_\epsilon\subset\{|t|\leq T_*\}$ be compact neighbourhoods of $A$ decreasing to $A$. By \eqref{weighted-duality}, $m(t)\Gamma((x,t);(y,s))=m(s)\widehat\Gamma((y,s);(x,t))$ for $s<t$. Tonelli's theorem gives
\[
 \int m\,V_{F_\epsilon}\,d\widehat\mu_A=\int m\,\widehat V_A\,d\mu_{F_\epsilon}.
\]
Since $\operatorname{supp}\widehat\mu_A\subset A\subset\operatorname{int}F_\epsilon$ and $V_{F_\epsilon}=1$ there, the left side equals $\widehat c_m(A)$. Since $\widehat V_A\leq1$, the right side is at most $c_m(F_\epsilon)$. Hence $\widehat c_m(A)\leq c_m(F_\epsilon)$. Because $\mu_F$ maximizes the weighted mass among admissible measures, $c_m(F)$ is the maximal weighted mass, and this quantity is continuous along decreasing compact sets: evaluation of the potentials at a point later than $F_\epsilon$ bounds the admissible masses, a weak limit of admissible measures is supported on the limit set, and it is admissible because the kernel is nonnegative and lower semicontinuous. Letting $\epsilon\downarrow0$ gives $\widehat c_m(A)\leq c_m(A)$. Interchanging the two operators gives the reverse inequality, and \eqref{weighted-capacity-duality} follows. The argument uses only \eqref{weighted-duality}, Tonelli's theorem and the balayage facts listed above, and therefore applies without change to the variable-coefficient operators in \cite{GII}.

Define
\[
 d\omega_A(\xi,\tau)=m_+^{-1}m(\tau)\,d\widehat\mu_A(\xi,\tau).
\]
Equation~\eqref{weighted-duality} gives
\[
 \Phi_A(x,t)=m_+^{-1}m(t)\widehat V_A(x,t)\leq1
                        \quad(|t|\leq T_*).
\]
Equations \eqref{weighted-extremum} and
\eqref{weighted-capacity-duality} give
$\omega_A(A)=m_+^{-1}c_m(A)\geq\chi^{-1}\operatorname{cap}(A)$.
This gives \eqref{dual-test}. We note that the argument only uses the two-sided comparison in \eqref{capacity-mass-comparison}, and not an equality between the capacity and $\mu_A(A)$.

\end{proof}

When $B=\overline B$, the estimate in the next lemma also follows from the behaviour of capacity under the dilations \eqref{sdl}. In the general case we use instead the small-time estimates for the Gramian.

\begin{lemma}[Capacity of a small box]\label{lem:box-capacity}
For the past boxes $C_R=\overline I_R\times[-R,0]$, there is $C$
such that
\begin{equation}\label{box-capacity}
 \operatorname{cap}(C_R)\leq CR^\nu
                         \quad(0<R\leq s_*/2).
\end{equation}
In particular $\operatorname{cap}(\{0\})=0$.
\end{lemma}

\begin{proof}
For $\zeta=(x,t)\in C_R$, the time difference from $(0,R)$ lies
in $[R,2R]$. The scaled Gramian bounds and
$|\Delta_{\sqrt R}^{-1}x|\leq1$ give
$\Gamma((0,R);\zeta)\geq cR^{-\nu}$.
For every admissible measure supported in $C_R$ this yields
$1\geq\Gamma\mu(0,R)\geq cR^{-\nu}\mu(C_R)$.
Taking the supremum proves the assertion.
\end{proof}

\subsubsection{The annular interaction estimate}\label{sec:interaction}

In this subsection we relate the capacity of an individual annulus to the mean decrements of one fixed potential. This estimate plays the role of the capacitary interaction step in the proofs in \cite{EG,GL88}. We note that in its proof we split the Riesz measure of the potential, so that no Dirichlet problem on a kernel level set is needed.

Fix $a>1$. Take $\alpha,\theta,C_H$ from
Theorem~\ref{thm:strong}. Choose a constant $c>1$ so that
\begin{equation}\label{choice-c}
 \theta c>1,\qquad b_+(\theta c)^{-\nu}<1,
\end{equation}
and put
\begin{equation}\label{interaction-scales}
 r=c\rho^{1/\nu},\qquad K=\frac{2(\alpha c)^\nu}{b_-},
 \qquad \sigma=K\rho,\qquad \kappa=(1-\log2)/2.
\end{equation}
We fix $\rho_0>0$ so small that $r\leq r_*$,
$\alpha r\leq s_*$, $2K\rho\leq\rho_*$, and $\Omega_{2K\rho}$ lies in
$-T_*<t<0$ whenever $0<\rho\leq\rho_0$.

\begin{lemma}[Capacity controlled by one mean decrement]
\label{lem:interaction}
Let $F\subset\{|t|\leq T_*\}$ be compact, and let
$V=V_F$, $\mu=\mu_F$, and $u=1-V$. Suppose that $A$ is compact,
$A\subset\operatorname{int}F$, and
\[
 A\subset\{\rho^{-1}\leq\Gamma_0\leq a\rho^{-1},\ t<0\},
 \qquad 0<\rho\leq\rho_0.
\]
Then, with $r$ and $K$ from \eqref{interaction-scales},
\begin{equation}\label{interaction}
 \frac{\operatorname{cap}(A)}{\rho}
       \frac1{|I_r|}\int_{I_r}u(x,-r)\,dx
 \leq C_I\bigl(\mathcal M_{K\rho}V-
                         \mathcal M_{2K\rho}V\bigr),
 \qquad C_I=\frac{\chi C_HK}{\kappa}.
\end{equation}
All constants are independent of $F,A,\rho$.
\end{lemma}

\begin{proof}
Split the Riesz measure into its restrictions to $\Omega_\sigma$ and $\Omega_\sigma^c$. Define
\begin{equation}\label{riesz-splitting}
 w=\Gamma(\mu|_{\Omega_\sigma}),\qquad
 v=1-\Gamma(\mu|_{\Omega_\sigma^c}).
\end{equation}
Then $u=v-w$ and $0\leq u\leq v\leq1$. Moreover $v$ is harmonic
in $\Omega_\sigma$, since its Riesz source there vanishes; by hypoellipticity it is smooth there. The choice of $\sigma$ ensures
\[
 \overline{\mathcal Q_r}\setminus\{0\}\subset\Omega_\sigma,
\]
including the earlier averaging slice. Indeed, applying
\eqref{slice-levels} at the spatial origin gives
\[
 \frac{(4\pi)^{-n/2}}{V(\alpha r)}
 \geq b_-(\alpha r)^{-\nu}=\frac2\sigma.
\]
Also $A\subset\mathcal E_{\theta r}$, because
\[
 \frac{(4\pi)^{-n/2}}{V(\theta r)}
 \leq b_+(\theta r)^{-\nu}
              =\frac{b_+(\theta c)^{-\nu}}\rho<\frac1\rho.
\]
Thus the strong Harnack theorem applies to $v$ and gives, on $A$,
\[
 v\geq C_H^{-1}\frac1{|I_r|}\int_{I_r}v(x,-r)\,dx
   \geq C_H^{-1}b,
 \qquad b=\frac1{|I_r|}\int_{I_r}u(x,-r)\,dx.
\]
Since $V=1$ in $\operatorname{int}F$, we have $u=0$ on $A$ and
hence $v=w$ there.

Apply Lemma~\ref{lem:weighted-capacity} to $A$. Positivity and
Tonelli's theorem, followed by \eqref{dual-test}, give
\begin{align*}
 \frac{b}{\chi C_H}\operatorname{cap}(A)
 &\leq\frac b{C_H}\omega_A(A)
 \leq\int_A v\,d\omega_A
 =\int_A w\,d\omega_A\\
 &=\int_{\Omega_\sigma}\left(\int_A\Gamma(\zeta;z)\,d\omega_A(\zeta)\right)
                                                  d\mu(z)
 \leq\mu(\Omega_\sigma).
\end{align*}
All the integrands in this computation are nonnegative.
Corollary~\ref{cor:mass-decrement} bounds this mass by
\[
 \mu(\Omega_\sigma)\leq\frac{\sigma}{\kappa}
       (\mathcal M_\sigma V-\mathcal M_{2\sigma}V).
\]
Substituting $\sigma=K\rho$ gives \eqref{interaction}.
\end{proof}

\section{Proof of Theorem~\ref{thm:wiener}}\label{sec:wiener}

The two implications in Theorem~\ref{thm:wiener} use capacity in different ways. If the Wiener
series converges, subadditivity shows that the balayage of the exterior
near the pole tends to zero, and the point is irregular. If the series
diverges, we select a separated subfamily of annuli, apply the
interaction estimate to one equilibrium potential, and use the strong
Harnack inequality to turn the resulting slice information into a
positive value at the pole. Finally, by comparing consecutive levels in a geometric sequence we obtain the equivalence with conditions (iii) and (iv).

\begin{proof}[Proof of Theorem~\ref{thm:wiener}]
The change of variables
$\xi=x-e^{(t_0-t)B}x_0$, $\tau=t-t_0$ preserves $\cL$, the kernel,
and Lebesgue measure, and hence also the capacity. We may therefore
take $z_0=0$. Fix the constants and small scales chosen in
Sections~\ref{sec:local-geometry}--\ref{sec:interaction}.

\smallskip\noindent
\emph{Convergence implies irregularity.}
Suppose the series in Theorem~\ref{thm:wiener}(ii) converges. For a large
integer $p$ and small $R$, let
\[
 F_R=U^c\cap C_R,\qquad
 H_{R,p}=F_R\cap(\{\Gamma_0\leq a^p\}\cup\{0\}).
\]
The kernel is set to zero for $t\geq0$ as before. The set $H_{R,p}$
is compact, and
$F_R\subset H_{R,p}\cup\bigcup_{j\geq p}A_j(0)$.
We note that an admissible measure has no atoms, since an atom would make its potential unbounded at the points immediately after the time of the atom. In particular, the pole carries no mass. The kernel
bounds on these sets, admissibility of their equilibrium measures,
and subadditivity give
\begin{align}
 V_{F_R}(0)
 &\leq V_{H_{R,p}}(0)+\sum_{j\geq p}V_{A_j(0)}(0)\notag\\
 &\leq a^p\operatorname{cap}(C_R)
          +a\sum_{j\geq p}a^j\operatorname{cap}(A_j(0)).
                                      \label{necessity-bound}
\end{align}
Letting first $R\downarrow0$, and using Lemma~\ref{lem:box-capacity} for the first term, and then $p\to\infty$, we find $\lim_{R\downarrow0}V_{F_R}(0)=0$.
By \eqref{balayage-regularity}, the point is irregular, which proves
the necessary implication by contraposition.

\smallskip\noindent
\emph{Divergence implies regularity.}
Suppose now that the series diverges. Choose an integer $q\geq2$
large enough that
\begin{equation}\label{sparse-choice}
 a^{1-q}<b_-c^{-\nu},\qquad a^q\geq2.
\end{equation}
At least one of its $q$ residue classes still gives a divergent
series. List its indices as $j_k=j_*+qk$, and put
\[
 \rho_k=a^{-j_k},\qquad r_k=c\rho_k^{1/\nu},
 \qquad\sigma_k=K\rho_k.
\]
Discard finitely many terms to satisfy all the small-scale
restrictions. For each $k$ with positive capacity choose a compact
$K_k\subset A_{j_k}(0)\setminus\{0\}$ such that
\begin{equation}\label{annulus-trimming}
 \operatorname{cap}(K_k)\geq\tfrac12\operatorname{cap}(A_{j_k}(0)).
\end{equation}
To obtain $K_k$, remove from $A_{j_k}(0)$ an open neighborhood of the
pole contained, in $t\leq0$, in a past
box of capacity at most $\tfrac12\operatorname{cap}(A_{j_k}(0))$.
Such boxes exist by \eqref{box-capacity}. The remainder is compact,
and subadditivity proves \eqref{annulus-trimming}. For a zero
capacity annulus take $K_k=\varnothing$.

Fix any small $R>0$. After discarding further finitely many terms,
the compact set
\begin{equation}\label{sparse-F}
 F=\{0\}\cup\bigcup_k K_k
\end{equation}
is contained in $U^c\cap C_R$.
Compactness and the inclusion follow from the shrinking assertion
in Lemma~\ref{lem:local-geometry}. The averaging slices
are disjoint from $F$. Indeed \eqref{slice-levels} and
\eqref{choice-c}--\eqref{sparse-choice} give on each closed slice
\begin{equation}\label{slice-gap}
 \frac{a^{1-q}}{\rho_k}
 <\Gamma_0(x,-r_k)<\frac1{\rho_k}
                 \quad(x\in\overline I_{r_k}).
\end{equation}
All earlier selected annuli have kernel at most
$a^{1-q}/\rho_k$, and all current or later ones have kernel at
least $1/\rho_k$. The pole has time zero, and hence none of these
sets meets the averaging slice.

Let $F_\epsilon$ be the closed Euclidean $\epsilon$-neighborhood of
$F$, for sufficiently small $\epsilon>0$, so that
$F_\epsilon\subset\{|t|\leq T_*\}$.
Set $V_\epsilon=V_{F_\epsilon}$ and $u_\epsilon=1-V_\epsilon$.
Each $K_k$ lies in $\operatorname{int}F_\epsilon$.
Lemma~\ref{lem:interaction} therefore yields
\begin{equation}\label{finite-interaction-sum}
 \sum_{k=k_0}^{N}\frac{\operatorname{cap}(K_k)}{\rho_k}
      \frac1{|I_{r_k}|}\int_{I_{r_k}}u_\epsilon(x,-r_k)\,dx
 \leq C_I\sum_{k=k_0}^{N}
      (\mathcal M_{\sigma_k}V_\epsilon-
                         \mathcal M_{2\sigma_k}V_\epsilon).
\end{equation}
Since $2\sigma_{k+1}\leq\sigma_k$ and the means are nonincreasing
in their parameter, the sum on the right is at most
\[
 \mathcal M_{\sigma_N}V_\epsilon-
             \mathcal M_{2\sigma_{k_0}}V_\epsilon
 \leq 1+2|\tr B|T_*.
\]
The last inequality is \eqref{potential-mean-bounds}, which holds also when $\tr B<0$ and the kernel $K_\rho$ in \eqref{adjusted-kernel} changes sign. We emphasize that the bound does not depend on $N$, $R$ or $\epsilon$.

For each fixed $k$, \eqref{slice-gap} places the averaging slice
outside $F$. Thus \eqref{compact-continuity} and bounded convergence
give
\[
 \frac1{|I_{r_k}|}\int_{I_{r_k}}u_\epsilon(x,-r_k)\,dx
 \longrightarrow b_k:=
 \frac1{|I_{r_k}|}\int_{I_{r_k}}(1-V_F)(x,-r_k)\,dx.
\]
First pass to the limit $\epsilon\downarrow0$ in the finite sum
\eqref{finite-interaction-sum}, then let $N\to\infty$. We obtain
\begin{equation}\label{weighted-average-sum}
 \sum_k\rho_k^{-1}\operatorname{cap}(K_k)b_k
              \leq C_I(1+2|\tr B|T_*)<\infty.
\end{equation}
By \eqref{annulus-trimming} the sum of the coefficients
$\rho_k^{-1}\operatorname{cap}(K_k)$ diverges. Since $b_k\geq0$,
\eqref{weighted-average-sum} implies $\liminf_k b_k=0$.
We note that only the slice averages have been passed to the limit, and that no convergence of $V_{F_\epsilon}(0)$ has been used.

We convert this information on the earlier slices into the
pointwise lower bound required by \eqref{balayage-regularity}. The semigroup
identity and $V_F=\Gamma\mu_F$ imply
\[
 V_F(0)\geq\int_{\R^n}p(0,x,r)V_F(x,-r)\,dx,
\]
because the right side retains the part of the source lying at
times less than $-r$. For $u=1-V_F$, therefore,
\[
 u(0)\leq\int_{\R^n}p(0,x,r)u(x,-r)\,dx
 \leq 1-\delta+
       C\frac1{|I_r|}\int_{I_r}u(x,-r)\,dx,
\]
where $0\leq u\leq1$ and \eqref{slice-mass} were used in the
last step. Taking $r=r_k$ along a subsequence with $b_k\to0$
gives $V_F(0)\geq\delta$. Since $F\subset U^c\cap C_R$,
\[
 V_{U^c\cap C_R}(0)\geq\delta.
\]
The constant $\delta$ is independent of the arbitrarily small
$R$. Equation~\eqref{balayage-regularity} proves regularity and
establishes the annular criterion in Theorem~\ref{thm:wiener}(ii).

\smallskip\noindent
\emph{Equivalence with the integral criteria.}
Monotonicity of $\operatorname{cap}(E_\rho(0))$ shows that the
integral in Theorem~\ref{thm:wiener}(iii) converges if and only if
the series
\[
 \sum_j a^j\operatorname{cap}(E_{a^{-j}}(0))
\]
converges.
Indeed, its integral over $[a^{-j-1},a^{-j}]$ lies between
\[
 (a-1)a^j\operatorname{cap}(E_{a^{-j-1}}(0))
 \quad\hbox{and}\quad
 (a-1)a^j\operatorname{cap}(E_{a^{-j}}(0)).
\]
Furthermore,
\[
 A_j(0)\subset E_{a^{-j}}(0)\subset\bigcup_{k\geq j}A_k(0),
\]
so monotonicity, subadditivity, and Tonelli give
\[
 \sum_{j\geq j_0}a^j\operatorname{cap}(A_j(0))
 \leq\sum_{j\geq j_0}a^j\operatorname{cap}(E_{a^{-j}}(0))
 \leq\frac a{a-1}
             \sum_{k\geq j_0}a^k\operatorname{cap}(A_k(0)).
\]
Together with the annular criterion already proved, this establishes
Theorem~\ref{thm:wiener}(iii).
For the fourth condition, substitute $\rho=V(s)$.
By \eqref{local-volume},
\begin{align*}
 \frac{d\rho}{\rho}&=(\nu+O(s))\frac{ds}{s},\\
 \frac{\operatorname{cap}(E_\rho(0))}{\rho^2}\,d\rho
 &=\frac{\operatorname{cap}(E_{V(s)}(0))}{V(s)}
      (\nu+O(s))\frac{ds}{s}.
\end{align*}
The multiplier $\nu+O(s)$ is positive and bounded above and below
for small $s$. This proves the equivalence of conditions (iii) and (iv).
Finally, every compact level set has finite capacity, as one sees by evaluating an admissible potential at a later point at which the kernel has a positive minimum on the set. Hence the integral over any interval of
parameters bounded away from zero is finite. Consequently, convergence
or divergence of the Wiener integral depends only on its behavior as
$\rho\downarrow0$.
\end{proof}

\medskip

\noindent \textbf{AI Declaration:} During the preparation of this manuscript, the author used ChatGPT
(OpenAI) as an editorial and computational aid. Its use was limited
to reorganizing portions of the manuscript, improving the exposition,
preparing LaTeX and diagrams, and checking selected algebraic
computations. The mathematical problems, ideas, results, and arguments
are those of the author, who assumes full responsibility for the manuscript.

\medskip

\noindent \textbf{Data Availability Statement:} Data sharing is not applicable to this article, as no datasets were generated or analyzed during the current study.

\end{document}